\documentclass[a4paper,11pt]{amsart}
\usepackage[english]{babel}
\usepackage[latin1]{inputenc}
\usepackage{mathrsfs}
\usepackage{amsfonts}
\usepackage{amssymb}
\usepackage{amsmath}
\usepackage{amsthm}
\usepackage{bbm}
\usepackage{vmargin}
\usepackage{epsfig}
\usepackage{fancyhdr}
\usepackage{extramarks}

\newcommand{\C}{\mathbf{C}}
\newcommand{\R}{\mathbf{R}}
\newcommand{\N}{\mathbf{N}}
\newcommand{\Z}{\mathbf{Z}}
\newcommand{\Q}{\mathbf{Q}}
\newcommand{\D}{\mathcal{D}}

\newcommand{\UH}{\mathbf{H}}
\newcommand{\Oint}{\mathcal{O}}

\newcommand{\GLR}{{{\rm GL}(2,\mathbf{R})}}

\newcommand{\Gal}{{\rm Gal}}
\newcommand{\SL}{{\rm SL}}
\newcommand{\GL}{{\rm GL}}
\newcommand{\sgn}{{\rm sgn}}

\newcommand{\PSLR}{{{\rm PSL}(2,\mathbf{R})}}
\newcommand{\PSL}{{\rm PSL}}

\newcommand{\cusp}{\mathcal{C}}
\newcommand{\eisen}{\mathcal{E}}
\newcommand{\ghn}{{\Gamma\backslash\mathbf{H}^n}}
\newcommand{\gh}{{\Gamma\backslash\mathbf{H}}}
\newcommand{\liea}{\mathfrak{a}}
\newcommand{\lieb}{\mathfrak{b}}

\newcommand{\liep}{\mathfrak{p}}
\newcommand{\idnorm}{\mathcal{N}}
\newcommand{\idphantom}{\phantom{{}=}}
\newcommand{\Tr}{{\rm Tr}}

\newcommand{\res}{{\rm res}}

\renewcommand{\Re}{{\rm Re }}
\renewcommand{\Im}{{\rm Im }}
\renewcommand{\phi}{\varphi}
\renewcommand{\epsilon}{\varepsilon}
\renewcommand{\newline}{\\ \strut}

\theoremstyle{plain}
\newtheorem{theorem}{Theorem}[section]
\newtheorem{lemma}[theorem]{Lemma}
\newtheorem{corollary}[theorem]{Corollary}

\newtheorem{proposition}[theorem]{Proposition}

\numberwithin{equation}{section}
\begin{document}

\title{Quantum Unique Ergodicity for Eisenstein Series on the Hilbert Modular Group over a Totally Real Field}
\thanks{The author was
  supported by a stipend (EliteForsk) from The Danish Agency for Science, Technology and Innovation}   
\author{Jimi L. Truelsen}
\address{Department of Mathematical Sciences, University of Aarhus, Ny Munkegade Building 1530, 8000 Aarhus C, Denmark}
\email{lee@imf.au.dk}

\begin{abstract}
W. Luo and P. Sarnak have proved the quantum unique ergodicity property for Eisenstein series on $\PSL(2,\Z) \backslash \UH$. Their result is quantitative in the sense that they find the precise asymptotics of the measure considered. We extend their result to Eisenstein series on $\PSL(2,\Oint) \backslash \UH^n$, where $\Oint$ is the ring of integers in a totally real field of degree $n$ over $\Q$ with narrow class number one, using the Eisenstein series considered by I. Efrat. We also give an expository treatment of the theory of Hecke operators on non-holomorphic Hilbert modular forms.
\end{abstract}

\maketitle

\section{Introduction}
Let $\UH$ denote the upper half-plane and $\Gamma$ be a Fuchsian group of the first kind. We equip the surface $\gh$ with the measure induced by the Poincar\'{e} metric $d\mu = \frac{dxdy}{y^2}$ on $\UH$. If $\Gamma$ is hyperbolic we know that the quotient $\gh$ is compact and that the Laplace-Beltrami operator $\Delta$ associated with this surface, given in local coordinates by $-y^2\left(\frac{\partial^2}{\partial x^2}+\frac{\partial^2}{\partial y^2}\right)$, has pure point spectrum
\begin{align*}
0 = \lambda_0 < \lambda_1\le\dots
\end{align*}
and that $\lambda_n \to \infty$ as $n \to \infty$. Inspired by quantum chaos (see \cite{sarnak95} and \cite{sarnak03} for excellent surveys) Z. Rudnick and P. Sarnak \cite{rudnick94} conjectured that
\begin{align}\label{questatement}
\vert \phi_j \vert^2 d\mu \to \frac{1}{\mu(\gh)}d\mu,
\end{align}
where $\{\phi_j\}$ is an orthonormal basis for $L^2(\gh)$ of eigenfunctions of $\Delta$ with $\Delta \phi_j = \lambda_j \phi_j$, and the convergence is in the weak-$*$ topology. This is known as the quantum unique ergodicity conjecture. It has been established by Y. Colin de Verdière \cite{colin85}, A. Shnirelman \cite{shnirelman74} and S. Zelditch \cite{zelditch87} that (\ref{questatement}) holds for a subsequence of full density.

If $\Gamma = \PSL(2,\Z)$ the quotient $\gh$ is no longer compact, and $\Delta$ does not have pure point spectrum. However, by the Weyl law it is known that 
\begin{align*}
\# \{j \in \N_0 \mid \vert t_j \vert \le T \} \sim
\frac{\mu(\Gamma\backslash \UH)}{4 \pi}T^2,
\end{align*}
where $\lambda_j = 1/4+t_j^2$ are the eigenvalues of $\Delta$. Thus the analogue of the quantum unique ergodicity conjecture is
\begin{align*}
\vert \phi_j \vert^2 d\mu \to \frac{3}{\pi}d\mu
\end{align*}
where $\{\phi_j\}$ is a complete set of orthonormal eigenfunctions of $\Delta$. It was proved in \cite{luo95} that if the $\phi_j$'s are Hecke eigenforms then the conjecture is true for a (large) subsequence of the full sequence and it has been proved by R. Holowinsky and K. Soundararajan \cite{holowinsky08} that the conjecture is implied by the Ramanujan-Petterson conjecture.

In \cite{luo95} a continuous spectrum analogue of the quantum unique ergodicity conjecture was proved. More precisely it was proved that for $A,B \subset \Gamma\backslash \UH$ compact and Jordan measurable, such that $\mu(B) \ne 0$, we have the limit
\begin{align*}
\frac{\int_A\vert E(z,1/2 + it)\vert^2 d\mu}{\int_B\vert E(z,1/2 + it)\vert^2 d\mu} \to \frac{\mu(A)}{\mu(B)}
\end{align*}
as $t \to \infty$, where $E(z,s)$ is the Eisenstein series on $\PSL(2,\Z)$. The authors even found explicit asymptotics for the measure $\vert E(z,1/2 + it)\vert^2 d\mu$ (in terms of integration of a continuous function with compact support). In this paper we generalize this result to Eisenstein series $E(z,s,m)$ (it will be defined in Section 11) on $\Gamma \backslash \UH^n$, where $\Gamma = \PSL(2,\Oint)$ and $\Oint$ is the ring of integers in a totally real field $K$ of degree $n$ over $\Q$ with narrow class number one. Note that instead of just one Eisenstein series as in the case of $\PSL(2,\Z)$ we have a family of Eisenstein series parametriced by $m \in \Z^{n-1}$.

We investigate the asymptotic behaviour of the measure $d\mu_{m,t} = \vert E(z,1/2+it,m)\vert^2d\mu$, where $\mu$ is the measure on $\ghn$ induced by the metric $\frac{dx_1\dots dx_n dy_1 \dots dy_n}{y_1^2 \dots y_n^2}$ on $\UH^n$:
\begin{theorem}\label{quelemma}
For $F \in C_c(\ghn)$ we have that
\begin{align*}
\frac{1}{\log t}\int_{\ghn}F(z)d\mu_{m,t}(z) \to \frac{\pi^{n} nR}{2 D \zeta_K(2)} \int_{\ghn}F(z)d\mu(z)
\end{align*}
as $t \to \infty$, where $\zeta_K$ denotes the Dedekind zeta-function and $D$ and $R$ denote the discriminant and regulator of $K$, respectively.
\end{theorem}
From this one easily deduces that:
\begin{theorem}\label{que}
Let $A,B \subset \ghn$ be compact and Jordan measurable, and assume
that $\mu(B) \ne 0$. Then
\begin{align*}
\frac{\mu_{m,t}(A)}{\mu_{m,t}(B)} \to \frac{\mu(A)}{\mu(B)}
\end{align*}
as $t \to \infty$.
\end{theorem}
To prove Theorem \ref{quelemma} we follow the same strategy as in \cite{luo95}. The idea in the proof is to find the asymptotics of $\int_{\ghn} f d\mu_{m,t}$, where $f$ is either an incomplete Eisenstein series or a Hecke eigenform, and then use the spectral decomposition of $L^2(\ghn)$. Estimates for various $L$-functions play a crucial role in the proof, and we will collect these results, as we go along. It should be mentioned that a similar result was shown in the case of a quadratic imaginary field with class number one in \cite{koyama00} using a subconvexity estimate (in the $t$-aspect) for the standard $L$-function proved in \cite{petridis01}.\newline

I would like to thank my advisor Morten S. Risager for suggesting this problem to me and for excellent guidance and supervision. I would also like to thank Akshay Venkatesh as well as the anonymous referee for useful comments. 

\section{Notation and Terminology}
Let $K$ be a totally real field of degree $n$ over $\Q$ and narrow
class number one (these are the standard assumptions which are usually
made to work with a non-adelic setup in textbooks such as
\cite{bump97} and \cite{garrett90}) and let $\Oint$ denote the ring of
integers in $K$. Here narrow class number one means that $\Oint$ is a principal ideal domain and that each non-zero ideal in $\Oint$ has a generator which is totally positive (this term is explained below).

Let
\begin{align}\label{galois}
\Gal(K/\Q) = \{\psi_1,\dots,\psi_n \}
\end{align}
with $\psi_1$ equal to the identity map on $K$. In this way we may regard $\Oint$ as a lattice in $\R^n$, by the
injection $\Oint \hookrightarrow \R^n$ defined by $a \mapsto
(a^{(1)},\dots,a^{(n)})$, where $a^{(j)} = \psi_j(a)$. Note that this embedding depends on the
choice of ordering of the elements in $\Gal(K/\Q)$ given in
(\ref{galois}).

We let $\Oint^\times$ denote the group of units in $\Oint$ and
$\Oint^* = \Oint - \{0\}$. The elements in $\Oint^*$ for which all the
embeddings are positive (such elements are called totally positive) will be denoted $\Oint_+$. We let $\Oint_+^\times = \Oint_+ \cap
\Oint^\times$ which clearly is a multiplicative group.

We let $\D$ denote the different, i.e. the inverse ideal of
\[
\D^{-1} = \{v \in K \mid \Tr(v \Oint) \subset \Z \}.
\]
It is a well known fact that $\D^{-1} \supset \Oint$ is a fractional ideal, and
since $K$ has narrow class number one there exists $\omega \in \Oint_+$ such
that $\D = (\omega) = \omega \Oint$ and $\D^{-1} = \omega^{-1} \Oint$.

It is well known that $\Oint$ is a free abelian group of rank $n$, and $\Oint^\times/\{\pm 1\}$ is a
free abelian group of rank $n-1$. In addition we know that for each $u \in \Oint^\times$ we have $\vert u^{(1)}\dots u^{(n)}\vert = 1$. We will assume that $\epsilon_1,\dots,\epsilon_{n-1} \in \R_+$ together with
$-1$ generate $\Oint^\times$. For later use let
\begin{align}\label{unitmatrix}
\begin{pmatrix}
e_{1,1}& \cdots & e_{1, n-1} & 1/n\\
\cdots & \cdots & \cdots & \cdots\\
e_{n,1}& \cdots & e_{n, n-1} & 1/n
\end{pmatrix}
=
\begin{pmatrix} \log \vert \epsilon_1^{(1)}\vert & \cdots & \log \vert \epsilon_1^{(n)}\vert \\
\cdots & \cdots & \cdots\\
\log \vert \epsilon_{n-1}^{(1)}\vert & \cdots & \log \vert \epsilon_{n-1}^{(n)}\vert \\
1 & \cdots & 1
\end{pmatrix}^{-1}.
\end{align}
Note that we have the relations
\begin{align}\label{unitrel1}
\sum_{j=1}^n e_{j,q} = 0,
\end{align}
and
\begin{align}\label{unitrel2}
\sum_{j=1}^n e_{j,q'}\log \vert \epsilon_q^{(j)} \vert = \delta_{q,q'}
\end{align}
for $q,q' = 1,\dots,n-1$.

We let $\UH$ denote the upper half plane of $\C$, i.e.
\begin{align*}
\UH = \{ z \in \C \mid \Im(z) > 0\}.
\end{align*}
We will use the convention $z = (z_1,\dots,z_n) \in \UH^n$ and
$z=(x,y)$ where $x = (x_1,\dots,x_n) \in \R^n$ and $y = (y_1,\dots,
y_n) \in \R_+^n$. Furthermore we will use the notation $dx= dx_1\dots
dx_n$ and $dy = dy_1\dots dy_n$.

We set $\Gamma = \PSL(2,\Oint) \subset \PSL(2,\R)$. This group is often referred to as the Hilbert modular group. The group $\Gamma$ does not in general imbed discretely in $\PSL(2,\R)$, but it does imbed discretely in $\PSL(2,\R)^n$ by the action on $\UH^n$ defined by
\begin{align*}
\pm\left( \begin{smallmatrix} a & b \\ c & d
    \end{smallmatrix}\right) z = \left(\frac{a^{(1)}z_1+b^{(1)}}{c^{(1)}z_1+d^{(1)}},\dots,\frac{a^{(n)}z_n+b^{(n)}}{c^{(n)}z_n+d^{(n)}}\right)
\end{align*}
which clearly is an extension of the classical action of $\PSL(2,\Z)$
on $\UH$ by M\"{o}bius transformations. For $\gamma = \pm\left( \begin{smallmatrix} a & b \\ c & d
    \end{smallmatrix}\right) \in \PSL(2,\Oint)$ we define $\gamma^{(j)} = \pm\left( \begin{smallmatrix} a^{(j)} & b^{(j)} \\ c^{(j)} & d^{(j)}
    \end{smallmatrix}\right)$.

If we regard $\UH^n$ as a Riemannian manifold with the metric
\begin{align*}
ds^2 = \frac{dx_1^2 + dy_1^2}{y_1^2}+\dots+\frac{dx_n^2 + dy_n^2}{y_n^2}
\end{align*}
the Laplace-Beltrami operator associated with this metric is
\begin{align*}
\Delta = \Delta_1 +\dots+ \Delta_n
\end{align*}
where $\Delta_j = -y_j^2\left(\frac{\partial^2}{\partial
    x_j^2}+\frac{\partial^2}{\partial y_j^2}\right)$. In the natural
    way the metric on $\UH^n$ transfers to the quotient $\ghn$. We also see that the $\Delta_j$'s induce symmetric and positive differential operators on
    $C_b^\infty(\ghn)$ which admit self-adjoint extensions (the
    Friedrichs extension). It is known that the quotient $\ghn$ has finite volume and as in the case $n = 1$ we will often regard functions on $\ghn$ as functions on the space $\UH^n$ which are invariant under $\Gamma$. The measure on $\ghn$ induced by the Riemannian metric is denoted $\mu$ and one can check that $d\mu = \frac{dxdy}{y_1^2 \dots y_n^2}$ in local coordinates.
\section{The Hecke $L$-function}
In the following it will be convenient to set $\rho_j(m) = \pi
\sum_{q=1}^{n-1}m_q e_{j,q}$ for $m \in \Z^{n-1}$. Let $\chi_m$ denote the following function on ${\C^*}^n$:
\begin{align}\label{expdef}
\chi_m(w)=\exp\left(i \pi \sum_{q=1}^{n-1}m_q \sum_{j=1}^n
  e_{j,q}\log\vert w_j \vert \right) = \prod_{j=1}^n\vert w_j\vert^{i \rho_j(m)}.
\end{align}
Clearly we can regard $\chi_m$ as a multiplicative function on $\Oint^*$ by the
usual embedding. For $\beta \in \Oint_+$ we note that $\chi_m(\beta)$ only depends on the
ideal $(\beta)$, so in this way we can regard $\chi_m$ as a multiplicative
function on the non-zero ideals in $\Oint$ (a so-called Gr\"{o}ssencharacter). Note also that for $m$ even, $\chi_m$ is trivial on $\Oint^\times$. We can now define the
Hecke $L$-function. It is defined by the series
\[
\zeta(s,m) = \sum_{\substack{\liea \subset \Oint\\ \liea \ne 0}} \frac{\chi_m(\liea)}{\idnorm(\liea)^s},
\]
which converges absolutely for $\Re(s) > 1$, and it can also be written as an Euler product over the prime ideals $\liep$, i.e.
\begin{align*}
\zeta(s,m) = \prod_{\liep}\left(1-\frac{\chi_m(\liep)}{\idnorm(\liep)^s} \right)^{-1}.
\end{align*}

The Hecke $L$-function has a meromorphic continuation to the entire complex plane. Furthermore $\zeta(s,m)$ is
entire if $m \ne 0$. The Dedekind zeta function
$\zeta(s,0)$ (sometimes also denoted $\zeta_K$) has a simple pole at $s=1$ with residue $\frac{2^{n-1}R}{\sqrt{D}}$ (cf. \cite{bump97} Section 1.7), and is holomorphic elsewhere. Here $D = \idnorm(\D) = \vert N(\omega) \vert$ is the discriminant
of $K$ and $R$ is the regulator of $K$, i.e. the absolute value of the determinant
\begin{align*}
\begin{vmatrix} \log \vert \epsilon_1^{(1)}\vert & \cdots & \log \vert \epsilon_1^{(n-1)}\vert \\
\cdots & \cdots & \cdots\\
\log \vert \epsilon_{n-1}^{(1)}\vert & \cdots & \log \vert \epsilon_{n-1}^{(n-1)}\vert \\
\end{vmatrix}.
\end{align*}

First we will make a convexity bound for the Hecke $L$-function on the
line $\Re(s) = \sigma$, where $\frac{1}{2}\le \sigma \le 1$. It is well known (see \cite{bump97} Theorem 1.7.2) that the Hecke
$L$-function $\zeta(s,m)$ satisfies the functional equation
\begin{align}\label{heckefunctional}
\xi(s,m) = \chi_m(\omega)i^{\Tr(\tau)}\xi(1-s,-m)
\end{align}
where $\xi(s,m)$ denotes the completed $L$-function defined by
\[
\xi(s,m) =
D^{s/2}\pi^{-ns/2}\zeta(s,m)\prod_{j=1}^n\Gamma\left(\frac{s+\tau_j-i\rho_j(m)}{2}\right),
\]
and $\tau = (\tau_1,\dots,\tau_n)$ is a binary vector depending on
$m$ with the property that
\begin{align}\label{taudef}
\chi_m((\beta)) = \chi_m(\beta) \prod_{j = 1}^n \sgn(\beta^{(j)})^{\tau_j}
\end{align}
for $\beta \in \Oint^*$.

Stirling's formula, i.e. the asymptotics of the $\Gamma$-function on vertical lines, plays a crucial role in the proof of Theorem \ref{quelemma}. For any $\sigma \in \R$ we have
\begin{align}\label{stirling1}
\Gamma(\sigma + it) \sim \sqrt{2 \pi}e^{-\pi \vert t\vert/2}\vert t \vert^{\sigma-1/2}
\end{align}
and
\begin{align}\label{stirling2}
\frac{\Gamma'(\sigma+it)}{\Gamma(\sigma+it)}
\sim \log\vert t \vert
\end{align}
as $\vert t \vert \to \infty$. Using the Phragm\'{e}n-Lindel\"{o}f principle (see \cite{iwaniec04}
Section 5.A), the functional equation (\ref{heckefunctional}) and Stirling's formula we easily derive the convexity bound
\begin{align}\label{convest}
\zeta(\sigma+it,m) \ll \vert t\vert^{\frac{n}{2}(1-\sigma)+\epsilon}
\end{align}
as $\vert t\vert \to \infty$, for any $\epsilon > 0$ and $\frac{1}{2}\le \sigma \le 1$.
Note that (\ref{convest}) gives the estimate
\begin{align*}
\zeta(1/2+it,m) \ll \vert t\vert^{\frac{n}{4}+\epsilon}
\end{align*}
for any $\epsilon > 0$. For later use it turns out that we need something slightly
better ($\frac{n}{4}-\epsilon$ in the exponent will do), i.e. we need a subconvexity estimate for $\zeta(s,m)$ on the critical
line. Such an estimate was proven by P. S\"{o}hne \cite{sohne97} (generalizing ideas due to D. R. Heath-Brown \cite{heathbrown78} and \cite{heathbrown88}):
\begin{theorem}\label{sohneest}
Let $\epsilon > 0$. Then
\begin{align*}
\zeta(1/2+it,m) \ll \vert t\vert^{\frac{n}{6}+\epsilon}
\end{align*}
as $\vert t\vert \to \infty$.
\end{theorem}
It is conjectured (and implied by the generalized Riemann hypothesis) that one in fact has
\begin{align*}
\zeta(1/2+it,m) \ll \vert t\vert^{\epsilon}
\end{align*}
for any $\epsilon > 0$ as $\vert t\vert \to \infty$.

It will also be necessary to estimate the logarithmic derivative of $\zeta(s,m)$ on the line $\Re(s) = 1$. We introduce a von Mangoldt type function on the non-zero ideals in
$\Oint$ defined by
\begin{align*}
\Lambda_m(\liea) =
\begin{cases}
\chi_m(\liea)\log \idnorm(\liep)  & \text{if } \liea = \liep^k\\
0 & \text{otherwise}
\end{cases},
\end{align*}
where $\liep$ denotes a prime ideal. For $\Re(s) > 1$ we see using the Euler product that
\begin{align*}
-\frac{\zeta'(s,m)}{\zeta(s,m)} &=
 -\sum_{\liep}\left(1-\frac{\chi_m(\liep)}{\idnorm(\liep)^s}\right)\frac{d}{ds}\left(\frac{1}{1-\frac{\chi_m(\liep)}{\idnorm(\liep)^s}}\right)\\
&= \sum_{\liep}\frac{\chi_m(\liep)\log
  \idnorm(\liep)}{\idnorm(\liep)^s\left(1-\frac{\chi_m(\liep)}{\idnorm(\liep)^s} \right)}\\
&=\sum_{\liep}\log \idnorm(\liep)
\sum_{k=1}^\infty\frac{\chi_m(\liep)^k}{\idnorm(\liep)^{sk}}\\
&=\sum_{\liea}\frac{\Lambda_m(\liea)}{\idnorm(\liea)^s}.
\end{align*}
Thus as in the case of the Riemann zeta-function the logarithmic derivative of $\zeta(s,m)$ can be written as a Dirichlet series.

To estimate the logarithmic derivative of the Hecke $L$-function we
need a zero-free region. By considering exponential sums one can obtain a zero-free region for
the Hecke $L$-function similar to Vinogradov's bound for the Riemann
zeta-function (see \cite{titchmarsh86} Chapter 6). This was done by M. Coleman \cite{coleman90}:
\begin{theorem}\label{coleman}
There exist positive constants $C$ and $L$ such that $\zeta(\sigma+it,m) \ne 0$ for
$\vert t \vert \ge L$ and $\sigma \ge 1 -\frac{C}{(\log \vert
  t\vert)^{2/3}(\log \log \vert t\vert)^{1/3}}$.
\end{theorem}
At present this is the best zero-free region we know, but the generalized Riemann hypothesis asserts that all zeros of $\zeta(s,m)$ in the critical strip $0 < \Re(s) < 1$ are on the line $\Re(s) = \frac{1}{2}$.

To obtain a sufficiently good estimate for the logarithmic derivative we follow
Landau's strategy (cf. \cite{titchmarsh86} Sections 3.9-3.11), which is based on the Borel-Carath\'{e}odory theorem. We remark that in order to use this approach it is necessary to estimate the Hecke $L$-function from below. To this end we consider the
following generalization of the M\"{o}bius function to non-zero ideals in
$\Oint$ defined by
\[
\mu(\liep_1^{\alpha_1}\dots\liep_k^{\alpha_k}) =
\begin{cases}
(-1)^k & \text{if } \alpha_1,\dots,\alpha_k \le 1\\
0 & \text{otherwise}
\end{cases}.
\]
The function $\mu$ has the following property (``M\"{o}bius inversion''):
\begin{equation}\label{mobiusinversion}
\sum_{\lieb \subset \liea}\mu(\liea) =
\begin{cases}
1 & \text{if } \lieb = \Oint\\
0 & \text{otherwise}
\end{cases},
\end{equation}
and the proof is the same as in the classical case (see
\cite{iwaniec04} Section 1.3). From this it is clear that
\begin{align}\label{invhl}
\frac{1}{\zeta(s,m)} = \sum_{\liea} \frac{\chi_m(\liea)\mu(\liea)}{\idnorm(\liea)^s}
\end{align}
for $\Re(s) > 1$. Thus
\begin{align*}
\frac{1}{\vert \zeta(\sigma + it,m) \vert} \le \zeta(\sigma,0)
\end{align*}
for $\sigma > 1$.

We have the following result due to Landau:
\begin{proposition}\label{landauthm}
Let $s = \sigma + it$ and assume that $\zeta(s,m) = O(e^{\phi(\vert t \vert)})$
for $\vert t \vert \ge L$ and $1-\theta(\vert t\vert)\le \sigma \le 2$ for some positive $L$, where $\phi(t)$ and $1/\theta(t)$ are positive increasing
functions defined for $t \ge L$ such that $\theta(t)\le 1$,
$\phi(t) \to \infty$ as $t \to \infty$ and $\phi(t)/\theta(t) =
o(e^{\phi(t)})$. Assume also that there exists a positive constant $C$
 such that $\zeta(s,m) \ne 0$ for $\vert t \vert \ge L$ and $\sigma \ge
1-C\frac{\theta(\vert t \vert)}{\phi(\vert t \vert)}$. Then
\begin{align*}
\frac{\zeta'(s,m)}{\zeta(s,m)} = O\left( \frac{\phi(\vert t \vert)}{\theta(\vert t \vert)}\right)
\end{align*}
and 
\begin{align*}
\frac{1}{\zeta(s,m)} = O\left(\frac{\phi(\vert t \vert)}{\theta(\vert t\vert)}\right)
\end{align*}
for $\vert t \vert \ge L+1$ and $\sigma \ge 1-\frac{C\theta(t)}{4\phi(t)}$.
\end{proposition}
Using Theorem \ref{coleman} we can apply Proposition \ref{landauthm} with $\phi(t) = (\log t)^{\frac{2}{3}}$ and $\theta(t) = (\log\log t)^{-\frac{1}{3}}$ to obtain the following:
\begin{corollary}\label{colemancor}
There exists a positive number $L$ such that for $\vert t \vert \ge L$ we have the
estimate
\begin{align*}
\frac{\zeta'(1+it,m)}{\zeta(1+it,m)} = O((\log t)^{\frac{2}{3}}(\log\log t)^{\frac{1}{3}}).
\end{align*}
\end{corollary}
In the same way we obtain an explicit lower bound for the Hecke $L$-function:
\begin{corollary}\label{convest2}
There exists a positive number $L$ such that for $\vert t \vert \ge L$ we have the
estimate
\begin{align*}
\frac{1}{\zeta(1+it,m)} = O((\log t)^{\frac{2}{3}}(\log\log t)^{\frac{1}{3}}).
\end{align*}
\end{corollary}

\section{Hecke Operators}
In this section we give an expository treatment of the theory of Hecke operators on non-holomorphic Hilbert modular forms. analogous to the treatment of Hecke operators on holomorphic Hilbert modular forms in \cite{bump97} Section 1.7 and \cite{garrett90}  Section 1.15.

We recall the abstract definition of the Hecke ring (see
\cite{shimura94}). We set $G = \GL(2,K)$, $\Gamma = \SL(2,\Oint)$ and let
$\mathfrak{D} \subset \GL(2,K)$ denote the $2 \times 2$ matrices with entries in $\Oint$ and totally positive determinant. The Hecke algebra
${{\rm R}(\Gamma,\mathfrak{D})}$ is the $\C$-vector space of finite formal sums $\sum_k
c_k \Gamma \alpha_k \Gamma$, where $\alpha_k \in \mathfrak{D}$ and $c_k \in
\C$. The addition in ${{\rm R}(\Gamma,\mathfrak{D})}$ is the obvious one,
while the multiplication is defined as follows. Let $\alpha, \beta \in
\mathfrak{D}$. It is well known that there exist distinct cosets
$\Gamma\alpha_1,\dots,\Gamma\alpha_d$ and
$\Gamma\beta_1,\dots,\Gamma\beta_{d'}$, where $\alpha_i, \beta_{i'} \in \mathfrak{D}$, such that $\Gamma \alpha \Gamma
= \cup_{i=1}^d \Gamma\alpha_i$ and $\Gamma \beta \Gamma =
\cup_{i'=1}^{d'} \Gamma\beta_{i'}$. We define $\Gamma \alpha \Gamma \cdot
\Gamma \beta \Gamma = \sum_{i,i'} \Gamma\alpha_i\beta_{i'}\Gamma$, which
clearly is independent of the choice of the $\alpha_i$'s and $\beta_{i'}$'s. We extend this multiplication in the obvious way, making ${{\rm R}(\Gamma,\mathfrak{D})}$ an algebra.

We can define a homomorphism from $\GLR_+$ to $\PSLR$ by mapping $\tau = \left( \begin{smallmatrix} a & b \\ c & d\end{smallmatrix}\right) \in \GLR_+$ to $w \mapsto \frac{aw+b}{cw+d}$ in $\PSLR$. Thus for $w \in \UH$ we simply define
\begin{align*}
\tau w = \frac{aw+b}{cw+d}.
\end{align*}
Therefore we get a natural map from $\mathfrak{D}$ to $\PSL(2,\R)^n$ and we see that ${{\rm R}(\Gamma,\mathfrak{D})}$ can be regarded as an algebra
of operators on $L^2(\ghn)$ (or even the vector space of automorphic functions) if we define $(\Gamma \alpha \Gamma f)(z) = \sum_{i=1}^d f(\alpha_i z)$.

Two double cosets $\Gamma \alpha \Gamma$ and $\Gamma \beta \Gamma$ are said to be equivalent if $\alpha = \eta \beta$ where $\eta = \left( \begin{smallmatrix} u & 0 \\ 0 & u\end{smallmatrix}\right)$ for some $u \in \Oint^\times$. Note that if $\alpha = \eta \beta$ then $\alpha^{(j)}z_j = \beta^{(j)}z_j$ for all $j=1,\dots,n$. Let $\nu \in \Oint_+$. Inspired by Hecke operators in the case of
holomorphic Hilbert modular forms (see \cite{garrett90}) we define
\begin{equation}\label{heckedef}
T_\nu f= \frac{1}{\sqrt{\vert N(\nu) \vert}}\sum_{\substack{\det
    \alpha = u \nu\\ u \in \Oint_+^\times}} \Gamma \alpha \Gamma f.
\end{equation} 
Here the sum should be taken over inequivalent double cosets.

We can use the class number one assumption to make this more explicit. Consider $\left( \begin{smallmatrix} a & b\\ c & d \end{smallmatrix}\right) \in \mathfrak{D}$. Write $a = r a'$ and $c = r c'$ where $a'$ and $c'$ are relative prime (i.e. $(a') + (c') = \Oint$). There exist $b', d' \in \Oint$ such that $a'd' - b'c' = 1$ and we see that
\begin{align*}
\begin{pmatrix}
d' & -b' \\
-c' & a'
\end{pmatrix}
\begin{pmatrix}
a & b \\
c & d
\end{pmatrix}
\end{align*}
is upper triangular. Thus for any $\alpha \in \mathfrak{D}$ we can find $\beta
\in \mathfrak{D}$, which is upper triangular and satisfies that $\Gamma  \alpha
= \Gamma  \beta$. Using this we can write the Hecke operator as follows
\begin{equation}\label{heckedef2}
T_\nu f(z) = \frac{1}{\sqrt{\vert N(\nu) \vert}}\sum_{\substack{ad = u \nu\\ u \in \Oint_+^\times}} \sum_{b \in \Oint/ (d)} f\left(\left( \begin{smallmatrix} a & b\\ 0 & d \end{smallmatrix}\right) z\right).
\end{equation} 
The outer sum is finite by unique factorization and the inner sum is
finite since $\Oint/(d)$ is a finite group. Thus $T_\nu f$ is well
defined.

If $\nu,\nu' \in \Oint_+$ and $\nu = u \nu'$ for some $u \in
\Oint_+^\times$ then by definition $T_\nu = T_{\nu'}$. Thus we define
$T_{(\nu)} = T_\nu$. Since we assumed that all ideals have a generator in
$\Oint_+$ there is a Hecke operator associated with each non-zero
ideal. Modifying Theorem 3.12.4 in \cite{goldfeld06} we obtain  that the Hecke operators are self-adjoint.

Now we will investigate the properties of the Hecke operators. We have the
following proposition:
\begin{proposition}\label{heckemult}
 Let $\nu_1, \nu_2 \in \Oint_+$ be relative prime. Then
\begin{align*}
T_{\nu_1}T_{\nu_2} = T_{\nu_1 \nu_2}.
\end{align*}
\end{proposition}
\proof Let $f \in L^2(\ghn)$. We have that
\begin{align*}
\sqrt{\vert N(\nu_1 \nu_2) \vert}&(T_{\nu_1}T_{\nu_2}f) (z)\\
&= \sum_{\substack{a_1 d_1 = u_1 \nu_1\\ u_1 \in \Oint_+^\times}}
\sum_{\substack{a_2d_2 = u_2 \nu_2\\ u_2 \in \Oint_+^\times}}
\sum_{\substack{b_1 \in \Oint/ (d_1)\\ b_2 \in \Oint/ (d_2)}} f\left( \left( \begin{smallmatrix} a_2 & b_2\\ 0 & d_2 \end{smallmatrix}\right)\left( \begin{smallmatrix} a_1 & b_1\\ 0 & d_1 \end{smallmatrix}\right)z\right)\\
&= \sum_{\substack{a_1d_1 = u_1 \nu_1\\ u_1 \in \Oint_+^\times}}
\sum_{\substack{a_2d_2 = u_2 \nu_2\\ u_2 \in \Oint_+^\times}}
\sum_{\substack{b_1 \in \Oint/ (d_1)\\ b_2 \in \Oint/ (d_2)}} f\left( \left( \begin{smallmatrix} a_1a_2 & b_1a_2+d_1b_2\\ 0 & d_1d_2 \end{smallmatrix}\right)z\right)\\
&= \sum_{\substack{ad = u \nu_1 \nu_2\\ u \in \Oint_+^\times}} \sum_{b \in \Oint/(d)} f\left( \left( \begin{smallmatrix} a & b\\ 0 & d \end{smallmatrix}\right)z\right)\\
&=\sqrt{\vert N(\nu_1 \nu_2) \vert}(T_{\nu_1 \nu_2}f)(z)
\end{align*}
where we have used the Chinese remainder theorem, i.e. that
\[
\Oint/((d_1)\cap (d_2)) \cong \Oint/(d_1)\oplus \Oint/(d_2)
\]
which holds since $(\nu_1) + (\nu_2) = \Oint$.\qed\newline

We need the following important lemma:
\begin{lemma}\label{heckeprim}
Let $p \in \Oint_+$ be a prime element. Then for any positive integers $k, k'$ we have
\[
T_{p^k}T_{p^{k'}} = \sum_{d=0}^{\min\{k,k'\}}T_{p^{k+k'-2d}}.
\]
\end{lemma}
\proof Let $f \in L^2(\Gamma\backslash \UH^n)$. We see that
\begin{align*}
\sqrt{\vert N(p^{k+k'}) \vert}(T_{p^k}T_{p^{k'}}f)(z) &= \sum_{\substack{l_1+l_2 = k\\ l_1'+l_2' = k'}}\sum_{\substack{b \in \Oint/ (p^{l_2})\\ b' \in \Oint/ (p^{l_2'})}}f\left( \left( \begin{smallmatrix} p^{l_1} & b\\ 0 & p^{l_2} \end{smallmatrix}\right)\left( \begin{smallmatrix} p^{l_1'} & b'\\ 0 & p^{l_2'} \end{smallmatrix}\right)z\right)\\
&= \sum_{\substack{l_1+l_2 = k\\ l_1'+l_2' = k'}} \sum_{\substack{b \in \Oint/ (p^{l_2})\\ b' \in \Oint/ (p^{l_2'})}} f\left( \left( \begin{smallmatrix} p^{l_1+l_1'} & bp^{l_2'}+b'p^{l_1}\\ 0 & p^{l_2+l_2'} \end{smallmatrix}\right)z\right).
\end{align*}
Removing common factors we get 
\[
\sum_{d=0}^{\min\{k,k'\}}\sum_{\substack{l_1+l_2 = k-d\\ l_1'+l_2' = k'-d\\ \min\{l_2',l_1\} = 0}} \sum_{b \in \Oint/ (p^{l_2})}\sum_{b' \in \Oint/ (p^{l_2'+d})}f\left( \left( \begin{smallmatrix} p^{l_1+l_1'} & b p^{l_2'}+b'p^{l_1}\\ 0 & p^{l_2+l_2'} \end{smallmatrix}\right)z\right).
\]
We note that as $(b,b')$ runs over all pairs in $\Oint/ (p^{l_2}) \times \Oint/ (p^{l_2'+d})$ the expression $b p^{l_2'}+b'p^{l_1}$ will assume each value in $\Oint/ (p^{l_2+l_2'})$ exactly $\vert N(p)\vert ^d$ times. Thus
\[
\sqrt{\vert N(p^{k+k'}) \vert}T_{p^k}T_{p^{k'}} = \sum_{d=0}^{\min\{k,k'\}}\vert N(p^d) \vert\sqrt{\vert N(p^{k+k'-2d})\vert}T_{p^{k+k'-2d}},
\]
and this proves the theorem.\qed\newline

Combining Proposition \ref{heckemult} and Lemma \ref{heckeprim} we get:
\begin{theorem}\label{heckethm}
Let $(\nu_1), (\nu_2)$ be non-zero ideals in $\Oint$. Then
\[
T_{(\nu_1)}T_{(\nu_2)} = \sum_{(d)\supset (\nu_1)+(\nu_2)}T_{(\nu_1)(\nu_2)/(d)^2}.
\]
In particular the Hecke operators commute.
\end{theorem}
From Lemma \ref{heckeprim} we obtain the following proposition:
\begin{proposition}\label{heckeexp}
Let $p \in \Oint_+$ be a prime element. Then for $k \in \N _0$ we have that
\begin{equation}\label{heckeeven}
T_{p^{2k}} = \sum_{l=0}^k (-1)^{k+l}\binom{k+l}{2l}{T_{p}}^{2l}
\end{equation}
and
\begin{equation}\label{heckeodd}
T_{p^{2k+1}} = \sum_{l=1}^{k+1} (-1)^{k+l+1}\binom{k+l}{2l-1}{T_{p}}^{2l-1}.
\end{equation}
\end{proposition}
\proof We first consider (\ref{heckeeven}). The claim certainly holds for $k = 0$ and $k = 1$. Now let $k' \ge 2$ be an integer and assume that the formula holds for $k \le k'$. Using Lemma \ref{heckeprim} we get
\begin{align*}
T_{p^{2k'+2}} &= T_{p^{2k'}}T_{p^2} - T_{p^{2k'}}-T_{p^{2k'-2}}\\
&= (T_p^2-2)T_{p^{2k'}}-T_{p^{2k'-2}}\\
&= T_p^{2k'+2}-(2k'+1)T_p^{2k'}+(-1)^{k'+1}+\\
&\idphantom\sum_{l=1}^{k'-1} (-1)^{k'+l+1}\left(2\binom{k'+l}{2l}+\binom{k'+l-1}{2l-2}-\binom{k'+l-1}{2l}\right){T_{p}}^{2l}\\
&= \sum_{l=0}^{k'+1} (-1)^{k'+l+1}\binom{k'+l+1}{2l}{T_{p}}^{2l}.
\end{align*}
By induction this proves (\ref{heckeeven}), and (\ref{heckeodd}) is proved by similar arguments.\qed

\section{The Fourier Expansion of an Automorphic Form}
An automorphic form $f$ is a formal eigenfunction of the Laplacians $\Delta_j$ (i.e. $f$ need not be square
integrable and we allow $f$ to be identically zero also) which satisfy the growth condition
\begin{align*}
f(z) = o(\exp(2 \pi y_j))
\end{align*}
as $y_j \to \infty$ for all $j =1,\dots,n$. This holds in particular if $f$ is square integrable. By construction we have $f(z+l) = f(z)$ for all $l \in \Oint$. Thus $f$ has a Fourier expansion (see \cite{huntley91}):
\begin{theorem}
Let $f$ be an automorphic form with Laplace eigenvalues $s_j(1-s_j)$. Then $f$ admits a Fourier expansion of the form
\begin{align*}
f(z) = \sum_{l \in \Oint} a_l(y)e(\Tr(\omega^{-1} l x)),
\end{align*}
where $e(x) = \exp(2\pi i x)$. Since $f(z)$ is an eigenfunction for the Laplacians $\Delta_1,\dots,\Delta_n$ the $l$-th Fourier coefficient $a_l(y)$ must satisfy the differential equations
\begin{equation}
\frac{\partial^2 a_l(y)}{\partial
  y_j^2}+\left(\frac{s_j(1-s_j)}{y_j^2}-4\pi^2\vert (\omega^{-1} l)^{(j)}\vert^2\right)a_l(y) = 0
\end{equation}
for $j=1,\dots,n$ and hence be of the form
\begin{equation*}
a_l(y) = c_l \sqrt{y_1\dots y_n}\prod_{j=1}^n
K_{s_j-\frac{1}{2}}(2\pi\vert (\omega^{-1} l)^{(j)}\vert y_j)
\end{equation*}
for $l \ne 0$. The zeroth Fourier coefficient can be written as a linear combination of $\prod_{j=1}^n y_j^{s_j}$ and $\prod_{j=1}^n y_j^{1-s_j}$. Furthermore the coefficients $c_l$ satisfy the bound
\begin{align*}
c_l \ll \exp(\epsilon \vert N(l)\vert)
\end{align*}
for any $\epsilon > 0$.
\end{theorem}
Here $K_\nu$ denotes the usual Macdonald Bessel function
\begin{align*}
K_\nu(y)=\frac{1}{2}\int_0^\infty\exp(-y(t+1/t)/2)t^{\nu-1}dt,
\end{align*}
which is defined for $y > 0$ and $\nu \in \C$. It is well known that
these functions decay exponentially as $y \to \infty$.

Note that if $f$ is automorphic with respect to $\Gamma$ then $f(z) = f(uz)$ for $u \in \Oint_+$, where
\begin{align*}
uz = (u^{(1)} z_1,\dots,u^{(n)}z_n),
\end{align*}
since all such $u$'s are squares of units (by the assumption that $K$ has narrow class number one). This implies that $c_{l} = c_{lu}$ for $l \in \Oint$ and $u \in \Oint_+^\times$.

A non-zero square integrable automorphic form $f$ is called a cusp form if
\begin{align}\label{cuspdef}
\int_F f(z)dx = 0.
\end{align}
Here
\begin{align*}
F = \{ t_1a_1+\dots +t_na_n \mid 0 \le t_j < 1\}
\end{align*}
where $a_1,\dots,a_n$ is a $\Z$-basis for $\Oint$ and each $a_j$ is regarded as a vector in $\R^n$ by the embedding $a_j \mapsto (a_j^{(1)},\dots,a_j^{(n)})$. We will refer to $F$ as the fundamental mesh for $\Oint$ and one can check that the definition of cuspidal is independent of the choice of $\Z$-basis. By the exponential decay of the Macdonald Bessel function one can deduce that $f$ must be of exponential decay as $y_j \to \infty$.

Using the Hilbert-Schmidt kernel from \cite{efrat87} Section II.9 one can prove using Lemma I.2.1 in \cite{efrat87} that the vector space of square integrable automorphic forms with given Laplace eigenvalues $\lambda_1,\dots,\lambda_n$ is finite dimensional (see \cite{huntley91} for bounds on the dimensions of the eigenspaces). Now define
\begin{align*}
\iota_j (z) = (z_1,\dots,z_{j-1},-\overline{z_j},z_{j+1},\dots,z_n)
\end{align*}
for $j= 1,\dots,n$. One easily checks that if $f$ is an automorphic form then so is $f \circ \iota_j$ with the same Laplace eigenvalues. Since the eigenspaces are finite dimensional this means that the eigenvalues of $\iota_j$ must be $\pm  1$. We also see that the Hecke operators, the $\Delta_j$'s and the $\iota_j$'s commute. Furthermore all these operators are self-adjoint. Hence we can choose a basis for the vector space spanned by cusp forms which consists of cusp forms that are also eigenfunctions for all the Hecke operators and all the $\iota_j$'s. These are called primitive cusp forms. Note that being an eigenfunction of the $\iota_j$'s is simply the same as saying that the function is either even or odd in each $x_j$.

\section{Hecke Eigenvalues and Automorphic Forms}
In this section we will study automorphic forms which are common eigenfunctions for all the Hecke operators. We first note that the identities derived in Theorem \ref{heckethm} and Proposition \ref{heckeexp} give similar identities for the 
Hecke eigenvalues:
\begin{theorem}\label{heckeeigen}
Assume that $f$ is a common eigenfunction for all the Hecke operators, i.e. that
\begin{align*}
T_{(\nu)} f= \lambda((\nu))f
\end{align*}
for all $\nu \in \Oint^*$. Then for $\nu_1,\nu_2 \in \Oint^*$ we have 
\begin{align}
\lambda((\nu_1))\lambda((\nu_2)) = \sum_{(d)\supset (\nu_1)+(\nu_2)}\lambda((\nu_1 \nu_2/d^2)).
\end{align}
For a prime element $p \in \Oint$ and $k \in \N _0$ we have that
\begin{equation}
\lambda((p^{2k})) = \sum_{l=0}^k (-1)^{k+l}\binom{k+l}{2l}{\lambda((p))}^{2l}
\end{equation}
and
\begin{equation}
\lambda((p^{2k+1})) = \sum_{l=1}^{k+1} (-1)^{k+l+1}\binom{k+l}{2l-1}\lambda((p))^{2l-1}.
\end{equation}
\end{theorem}
Using the identities above, we can derive a connection between the Fourier coefficients of $T_{(\nu)}f$ and $f$, where $f$ is a primitive cusp form:
\begin{theorem}
Let $f$ be a primitive cusp form with Laplace eigenvalues $s_j(1-s_j)$, and assume that $f$ has the Fourier expansion
\[
f(z) = \sum_{l \in \Oint^*}c_l\sqrt{y_1\dots
 y_n}\left(\prod_{j=1}^nK_{s_j- \frac{1}{2}}(2\pi \vert (\omega^{-1}
 l)^{(j)}\vert y_j) \right)e(\Tr(\omega^{-1} l x)).
\]
Then the $l$-th Fourier coefficient of $T_{(\nu)}f$ is
\[
\sum_{\substack{d\mid \gcd(l',\nu)\\ l'\nu = d^2l}}c_{l'}
\]
for $\nu \in \Oint_+$. In particular $c_{\nu u} = \lambda((\nu))c_u$ for $u
\in \Oint^\times$.
\end{theorem}
\proof We apply $T_{\nu}$ on the Fourier expansion
\begin{align*}
&\sqrt{\vert N(\nu) \vert}T_\nu f(z) = \sum_{l' \in \Oint^*}c_{l'}\sum_{\substack{ad = u \nu\\ u \in \Oint_+^\times}}\sqrt{\frac{\vert a^{(1)} \vert}{\vert d^{(1)} \vert}y_1\dots \frac{\vert a^{(n)} \vert}{\vert d^{(n)} \vert} y_n} \times\\
&\phantom{==}\left(\prod_{j=1}^nK_{s_j- \frac{1}{2}}(2\pi \vert (\omega^{-1}l'a/d)^{(j)}\vert y_j) \right)\sum_{b \in \Oint/ (d)} e\left(\Tr \left(\omega^{-1} l' (ax+b)/d\right)\right),
\end{align*}
where by abuse of notation
\begin{align*}
(ax+b)/d = ((a^{(1)}x_1+b^{(1)})/d^{(1)},\dots, (a^{(n)}x_n+b^{(n)})/d^{(n)}).
\end{align*}
Now if $d \mid l'$ then
\begin{align*}
\sum_{b \in \Oint/ (d)} e\left(\Tr \left(\omega^{-1} l' \frac{b}{d} \right)\right) = \vert N(d) \vert.
\end{align*}
If $d \nmid l'$ there exist $b' \in \Oint/(d)$ such that $\Tr \left(\omega^{-1} l' \frac{b'}{d} \right) \notin \Z$. Thus
\begin{align*}
\sum_{b \in \Oint/ (d)} e\left(\Tr\left(\omega^{-1} l' \frac{b}{d} \right)\right) &= \sum_{b \in \Oint/ (d)} e\left(\Tr\left(\omega^{-1} l' \frac{b+b'}{d} \right)\right) \\
&= e\left(\Tr\left(\omega^{-1} l' \frac{b'}{d} \right)\right)\times\\
&\phantom{{}=}\sum_{b \in \Oint/ (d)} e\left(\Tr \left(\omega^{-1} l' \frac{b}{d} \right)\right).
\end{align*}
But this implies that
\begin{align*}
\sum_{b \in \Oint/ (d)} e\left(\Tr\left(\omega^{-1} l' \frac{b}{d} \right)\right) = 0.
\end{align*}
Thus
\begin{align*}
  \sqrt{\vert N(\nu) \vert}T_{(\nu)}f(z) &= \sum_{l' \in
    \Oint^*}c_{l'} \sum_{\substack{ad = u\nu\\ u \in
      \Oint_+^\times}}\sqrt{\frac{\vert a^{(1)} \vert}{\vert d^{(1)}
      \vert}y_1\dots \frac{\vert a^{(n)} \vert}{\vert d^{(n)} \vert}
    y_n} \times \\
 &\phantom{=} \left(\prod_{j=1}^n K_{s_j-
      \frac{1}{2}} (2\pi \vert (\omega^{-1} l'\nu/d^2)^{(j)}\vert y_j)
    \right)\times\\
&\phantom{=}\sum_{b \in \Oint/ (d)} e\left(\Tr\left(\omega^{-1} l'
      (ax+b)/d\right)\right)\\
  &= \sum_{l' \in \Oint^*}c_{l'}\sum_{d \mid \gcd(l',\nu)}\vert N(d)
  \vert\sqrt{\frac{\vert \nu^{(1)} \vert}{\vert d^{(1)} \vert^2}y_1\dots
    \frac{\vert \nu^{(n)} \vert}{\vert d^{(n)} \vert^2}
    y_n}\times\\
&\phantom{=} \left(\prod_{j=1}^nK_{s_j- \frac{1}{2}}(2\pi \vert
    (\omega^{-1} l'\nu/d^2)^{(j)}\vert y_j) \right)\times\\
&\phantom{=} e\left(\Tr\left(\omega^{-1} l'
      \frac{\nu}{d^2}x\right)\right).
\end{align*}
From this it is clear that the $l$-th Fourier coefficient is
\[
\sum_{\substack{d\mid \gcd(l',\nu)\\ l'\nu = d^2l}}c_{l'}.
\]
\qed

\section{The Fundamental Domain for $\Gamma_\infty$}
Before we can prove the functional equation for the standard $L$-function we need a fundamental domain for $\Oint_+^\times\backslash \R_+^n$ and this immediately gives us a fundamental domain for $\Gamma_\infty$ as well.

Let $F$ denote the interior of the fundamental mesh of the lattice $\Oint$ in $\R^n$ given by the embedding defined earlier. Let $\Gamma_\infty$ denote the stabilizer subgroup at $\infty$, i.e.
\[
\Gamma_\infty = \left\{ \pm\left( \begin{smallmatrix} u & l \\ 0 & u^{-1}
    \end{smallmatrix}\right) \mid u \in \Oint^\times,l \in \Oint \right\}.
\]
From \cite{siegel65} we know the fundamental domain for
$\Gamma_\infty$:
\begin{proposition}\label{fundprop}
The set
\[
F_\infty=\left\{ z \in \UH^n \mid x \in F,y \in U_\infty\right\},
\]
is a fundamental domain for $\Gamma_\infty$. Here $U_\infty \subset
\R_+^n$ is a fundamental domain for $\Oint_+^\times \backslash \R_+^n$. Explicitly we can choose $U_\infty$ to be the preimage of
\begin{align*}
\R_+ \times [-1,1]^{n-1} \subset \R_+ \times \R^{n-1}
\end{align*}
under the map (defined on $\R_+^n$)
\begin{align*}
y \mapsto \left( \prod_{j=1}^n y_j,
  \sum_{j=2}^n (e_{j,1}-e_{1,1}) \log \frac{y_j}{\sqrt[n]{\prod_{i=1}^n
  y_i}},\dots,\sum_{j=2}^n (e_{j,n-1}-e_{1,n-1}) \log \frac{y_j}{\sqrt[n]{\prod_{i=1}^n
  y_i}} \right),
\end{align*}
which is injective.
\end{proposition}
Let $\widetilde{y}$ denote the image of $y$ under the map above. Note that we have the relations
\begin{align}\label{fundequa}
\sum_{j=2}^n \widetilde{y}_j \log \vert \epsilon_{j-1}^{(k)} \vert = \log \frac{y_k}{\sqrt[n]{\widetilde{y}_1}}
\end{align}
for $k = 2,\dots,n$ which follows since
\[
\begin{pmatrix}
e_{2,1}& \cdots & e_{2, n-1} \\
\cdots & \cdots & \cdots \\
e_{n,1}& \cdots & e_{n, n-1} 
\end{pmatrix}^{-1}
=
\begin{pmatrix} \log \vert \epsilon_1^{(2)}\vert & \cdots & \log \vert \epsilon_1^{(n)}\vert \\
\cdots & \cdots & \cdots\\
\log \vert \epsilon_{n-1}^{(2)}\vert & \cdots & \log \vert \epsilon_{n-1}^{(n)}\vert
\end{pmatrix}(I_{n-1}+E_{n-1}).
\]
Here $I_{n-1}$ denotes the $(n-1) \times (n-1)$ identity matrix and
$E_{n-1}$ is the $(n-1) \times (n-1)$ matrix with all entries equal to
$1$. Inserting (\ref{fundequa}) in (\ref{expdef}) we get the relation
\begin{align}\label{ytilderel}
\chi_m (y) = \exp\left( i\pi \sum_{q=1}^{n-1} m_q \widetilde{y}_{q+1}\right).
\end{align}
Note also that by (\ref{fundequa}) the ratios $y_j/y_i$ are bounded.

Later we want to integrate so-called incomplete Eisenstein
series. To do so it will be convenient to use the transformation from
Proposition \ref{fundprop} and for that purpose we need to know the
Jacobian determinant:
\begin{lemma}\label{jacob}
The numerical value of the Jacobian determinant of the map in
Proposition \ref{fundprop} is $R^{-1}$ where $R$ is the
regulator of $K$.
\end{lemma}
\proof Let $\Omega$ denote the Jacobian matrix. Note that
\begin{align*}
\frac{\partial \widetilde{y}_1}{\partial y_j} = \frac{\widetilde{y}_1}{y_j}
\end{align*}
and
\begin{align*}
\frac{\partial \widetilde{y}_{k+1}}{\partial y_j} = \frac{1}{y_j}\sum_{j' = 2}^n \delta_{j,j'}(e_{j',k}-e_{1,k})-\frac{1}{n y_j}\sum_{j' = 2}^n(e_{j',k}-e_{1,k})
\end{align*}
for $k = 1,\dots,n-1$. Thus the $y_j$'s cancel in the Jacobian determinant and we get
\begin{align*}
\det(\Omega) = \begin{vmatrix} 1 & 1 & \cdots & 1\\
A_1 & e_{2,1}-e_{1,1} + A_1 & \cdots & e_{n,1}-e_{1,1} + A_1\\
\cdots & \cdots & \cdots & \cdots\\
A_{n-1} & e_{2,n-1}-e_{1,n-1} + A_{n-1} & \cdots & e_{n,n-1}-e_{1,n-1} + A_{n-1}
\end{vmatrix}
\end{align*}
where $A_k = -\frac{1}{n}\sum_{j=2}^n(e_{j,k}-e_{1,k})$ for $1 \le k \le
n-1$. By recursively subtracting column $j-1$ from column $j$ we do not
change the determinant. Expanding by minors in the first row (which
has $1$ in the first entry and $0$ in the other entries) we see that
\begin{align*}
\det(\Omega)= \det([e_{j+1,k}-e_{j,k}]_{1\le j,k \le n-1}).
\end{align*}
Now using
a similar trick on the matrix 
\begin{align*}
\begin{pmatrix}
e_{1,1}& \cdots & e_{1, n-1} & 1/n\\
\cdots & \cdots & \cdots & \cdots\\
e_{n,1}& \cdots & e_{n, n-1} & 1/n
\end{pmatrix}
=
\begin{pmatrix} \log \vert \epsilon_1^{(1)}\vert & \cdots & \log \vert \epsilon_1^{(n)}\vert \\
\cdots & \cdots & \cdots\\
\log \vert \epsilon_{n-1}^{(1)}\vert & \cdots & \log \vert \epsilon_{n-1}^{(n)}\vert \\
1 & \cdots & 1
\end{pmatrix}^{-1}
\end{align*}
recursively subtracting row $k-1$ from row $k$ we see that
\begin{align*}
\pm\frac{\det(\Omega)}{n} = \begin{vmatrix} \log \vert \epsilon_1^{(1)}\vert & \cdots & \log \vert \epsilon_1^{(n)}\vert \\
\cdots & \cdots & \cdots\\
\log \vert \epsilon_{n-1}^{(1)}\vert & \cdots & \log \vert \epsilon_{n-1}^{(n)}\vert \\
1 & \cdots & 1
\end{vmatrix}^{-1}.
\end{align*}
But the determinant on the right-hand side is $\pm nR$ (see e.g. \cite{mazanti02}).\qed

\section{The Standard $L$-function}
In this section we will consider the $L$-function associated with a
primitive cusp form -- the so-called standard $L$-function -- and show that it has a functional equation.

For a primitive cusp form $\phi$ with Hecke eigenvalues $\lambda(\liea)$ we
consider the $L$-function (defined for $\Re(s) > \frac{3}{2}$)
\[
L(s,\phi,m) = \sum_{\liea \ne 0}\frac{\chi_m(\liea)\lambda(\liea)}{\idnorm(\liea)^s}.
\]
It should be mentioned that one often uses the notation $L(s,\phi \otimes \chi_m)$ instead of $L(s,\phi, m)$.

If we use the relations from Theorem \ref{heckeeigen} we can write
$L(s,\phi,m)$ as the Euler product
\[
L(s,\phi, m) = \prod_{\liep}\frac{1}{1-\frac{\chi_{m}(\liep)\lambda(\liep)}{\idnorm(\liep)^s}+\frac{\chi_{m}(\liep)^2}{\idnorm(\liep)^{2s}}}
\]
where the product is taken over all prime ideals.

Before we go on we need the following result:
\begin{lemma}\label{vanishlemma}
Let $f$ be a formal eigenfunction of the Laplacians $\Delta_1,\dots,\Delta_n$ with eigenvalues
$\lambda_1,\dots,\lambda_n$. Assume that $f(iy) = 0$ for all
$y \in \R_+^n$ where $iy = (iy_1,\dots,iy_n)$. Assume also that
\begin{align}\label{oddcond}
\frac{\partial f}{\partial x_j}(z_1,\dots,z_{j-1},iy_j,z_j,\dots,z_n) = 0
\end{align}
for all $z_{j'} \in \UH$ with $j' \ne j$, $y_j \in \R_+$ and $j
= 1,\dots,n$. Then $f(z) = 0$ for all $z \in \UH^n$.
\end{lemma}
\proof
Since $f$ is an eigenfunction of the $\Delta_j$'s which are elliptic
differential operators we conclude that $f$ must be real analytic. Hence it
suffices to prove that
\begin{align*}
\frac{\partial^{\vert a + b \vert} f}{\partial x_1^{a_1}\dots \partial
  x_n^{a_n}\partial y_1^{b_1}\dots \partial
  y_n^{b_n}}(iy) = 0
\end{align*}
for all $a = (a_1,\dots,a_n), b = (b_1,\dots,b_n)  \in \N_0^n$ and $y \in
\R_+^n$ -- note that we use the notation $\vert a \vert = \sum_{j=1}^n a_j$. But clearly this would follow if we could prove that
\begin{align*}
\frac{\partial^{\vert a \vert} f}{\partial x_1^{a_1}\dots \partial
  x_n^{a_n}}(iy) = 0
\end{align*}
for all $(a_1,\dots,a_n) \in \N_0^n$ and $y \in
\R_+^n$.

If $a_j \in \{0,1\}$ for some $j$, the result follows immediately from
(\ref{oddcond}). Now assume that the result holds if for some $j$ we
have $a_j \le q$, $q \ge 2$. Consider $(a_1,\dots,a_n) \in \N_0^n$
such that $\min\{a_1,\dots,a_n\} = q + 1$; say $a_1 = q+1$. Then we see that
\begin{align*}
\frac{\partial^{\vert a \vert} f(iy)}{\partial x_1^{a_1}\dots \partial
  x_n^{a_n}} &= \frac{\partial^{\vert a \vert-2} }{\partial x_1^{a_1-2}\partial x_2^{a_2}\dots \partial
  x_n^{a_n}}\left(-\frac{1}{y_1^2}\Delta_1 f(iy)-\frac{\partial^2
  f(iy)}{\partial y_1^2}\right)\\
&=-\frac{\lambda_1}{y_1^2}\frac{\partial^{\vert a \vert-2} f(iy)}{\partial x_1^{a_1-2}\partial x_2^{a_2}\dots \partial
  x_n^{a_n}}-\frac{\partial^2
  }{\partial y_1^2}\frac{\partial^{\vert a \vert-2}f(iy)}{\partial x_1^{a_1-2}\partial x_2^{a_2}\dots \partial
  x_n^{a_n}}\\
&=0
\end{align*}
by induction. This proves the lemma.\qed\newline

Now we can extend the holomorphic function $L(s,\phi,m)$ to an entire
function with a functional equation of the usual form:
\begin{theorem}\label{slffe}
Let $\phi$ be a primitive cusp form with Laplace eigenvalues $\frac{1}{4} + r_j^2$
and Hecke eigenvalues $\lambda(\liea)$. Then $L(s,\phi,m)$ has an analytic continuation to the entire
complex plane and it satisfies the functional equation 
\begin{align}\label{slffunctional}
\Lambda(s,\phi,m) = (-1)^{\Tr(\kappa)}\chi_{2m}(\D)\Lambda(1-s,\phi,-m)
\end{align}
where 
\begin{align*}
&\Lambda(s,\phi,m) =
 D^s \pi^{-ns}L(s,\phi,m)\times\\
&\phantom{===}\prod_{j=1}^n\Gamma\left(\frac{s+\kappa_j+ir_j-i\rho_j(m)}{2}\right)\Gamma\left(\frac{s+\kappa_j-ir_j-i\rho_j(m)}{2}\right)
\end{align*}
and $\kappa_j = 0$ if $\phi$ is even in $x_j$ and $\kappa_j = 1$ if $\phi$ is odd in $x_j$. 
\end{theorem}
\proof Consider the function
\begin{align*}
f(z) &=\frac{1}{(2 \pi i)^{\Tr(\kappa)}} \frac{\partial^{\Tr(\kappa)}\phi}{\partial x_1^{\kappa_1}\dots\partial
  x_n^{\kappa_n}}(z)\\
&=\sum_{l \in \Oint^*} c_l e(\Tr(l x/\omega))
\prod_{j=1}^n   \left(\frac{l^{\kappa_j}}{\omega^{\kappa_j}}\right)^{(j)}\sqrt{y_j}K_{ir_j}(2\pi \vert (l/\omega)^{(j)}\vert
  y_j),
\end{align*}
which is even in all the $x_j$-variables. For $\Re(s)$ large (this ensures that we can use the Fourier expansion) consider the integral
\begin{align*}
&\frac{\chi_{m}(\D)}{D^s}\int_{\Oint_+^\times \backslash \R_+^n} f(iy) \prod_{j=1}^n y_j^{s -
  i \rho_j(m) + \kappa_j -3/2}dy\\
&\phantom{0000}= \sum_{\liea \subset \Oint}
  \frac{\chi_{m}(\liea)\lambda(\liea)}{\idnorm(\liea)^s}\prod_{j=1}^n\int_0^\infty
  K_{ir_j}(2 \pi  y_j) y_j^{s-i\rho_j(m)+ \kappa_j-1}dy_j\times\\
&\phantom{0000=}\sum_{\beta
  \in \Oint_+^\times \backslash \Oint^\times} c_\beta \prod_{j=1}^n(\sgn(\beta^{(j)}))^{\tau_j}\displaybreak\\
&\phantom{0000}=
  L(s,\phi,m)\prod_{j=1}^n\frac{\Gamma\left(\frac{s+\kappa_j+ir_j-i\rho_j(m)}{2}\right)\Gamma\left(\frac{s+\kappa_j-ir_j-i\rho_j(m)}{2}\right)}{4\pi^{s+\kappa_j-i\rho_j(m)}} \times\\
&\phantom{0000=}\sum_{\beta
  \in \Oint_+^\times \backslash \Oint^\times} c_\beta \prod_{j=1}^n(\sgn(\beta^{(j)}))^{\tau_j}
\end{align*}
where $\tau$ is the binary vector satisfying (\ref{taudef}). Note that we have used the formula (see \cite{bump97} Lemma 1.9.1)
\begin{align*}
\int_0^\infty K_\nu(2\pi y) y^{s-1}dy = \frac{\Gamma\left(\frac{s+\nu}{2}\right)\Gamma\left(\frac{s-\nu}{2}\right)}{4\pi^s}
\end{align*}
which is valid for $\Re(s) > \vert \Re(\nu) \vert$. That the integral above is convergent follows from the fact that $f(iy) =
  \frac{(-1)^{\Tr(\kappa)}}{\prod_{j=1}^n y_j^{2\kappa_j}}f(i/y)$ (we use the notation $1/y = (1/y_1,\dots,1/y_n)$).

If we can prove that
\begin{align}\label{fcnonzero}
\sum_{\beta
  \in \Oint_+^\times \backslash \Oint^\times} c_\beta \prod_{j=1}^n(\sgn(\beta^{(j)}))^{\tau_j} \neq 0
\end{align}
we have the analytic continuation since
\begin{align*}
\int_{\Oint_+^\times \backslash \R_+^n} f(iy) \prod_{j=1}^n y_j^{s -
  i \rho_j(m) +\kappa_j -3/2}dy
\end{align*}
 is an entire function in $s$ (due to exponential decay of $f$ in
 the $y_j$-variables). So let us assume that 
\begin{align}
\sum_{\beta
 \in \Oint_+^\times \backslash \Oint^\times} c_\beta \prod_{j=1}^n(\sgn(\beta^{(j)}))^{\tau_j} = 0.
\end{align}
This implies that the integral considered above vanishes for all $s$ and $m \in \Z^{n-1}$. But using the structure of $U_ \infty$ we see that ($\widetilde{f}$ is $y \mapsto f(iy)$ composed with the inverse of the map in Proposition \ref{fundprop})
\begin{align*}
\int_{\Oint_+^\times \backslash \R_+^n} f(iy) &\prod_{j=1}^n y_j^{s -
  i \rho_j(m) + \kappa_j -3/2}dy =\\
&R\int_{-1}^1\dots \int_{-1}^1\int_0^\infty \widetilde{f}(\widetilde{y})\widetilde{y}_1^{s-3/2+\Tr(\kappa)/n}\times\\
&\exp\left(\sum_{q=1}^{n-1}\sum_{j=2}^n(\kappa_{q+1}-\kappa_1)\widetilde{y}_j\log \vert \epsilon_{j-1}^{(q+1)}\vert\right)\exp\left(- i\pi \sum_{q=1}^{n-1} m_q \widetilde{y}_{q+1}\right)d\widetilde{y},
\end{align*}
where we have used (\ref{ytilderel}). Since this holds for all $m$ we must have $f(iy) = 0$ for all
$y \in \R_+^n$. We also have that $f$ is a formal eigenfunction of
the $\Delta_j$'s and since $f$ is even in all the $x_j$-variables condition
(\ref{oddcond}) in Lemma \ref{vanishlemma} is also satisfied. Thus we
conclude that $f$ is identically $0$. But by the Fourier expansion of
$f$ this implies that $c_l = 0$ for all $l \in \Oint^*$ which contradicts
that $\phi$ is a primitive cusp form and hence non-zero.

Now we prove the functional equation. As remarked earlier $f(iy) =
  \frac{(-1)^{\Tr(\kappa)}}{\prod_{j=1}^n y_j^{2\kappa_j}}f(i/y)$. From this one
  easily deduces that 
\begin{align*}
\int_{\Oint_+^\times \backslash \R_+^n} f(iy) &\prod_{j=1}^n y_j^{s -
  i \rho_j(m) + \kappa_j -3/2}dy\\
&= (-1)^{\Tr(\kappa)}\int_{\Oint_+^\times \backslash \R_+^n} f(i/y) \prod_{j=1}^n y_j^{
s -  i \rho_j(m) - \kappa_j - 3/2}dy\\
&= (-1)^{\Tr(\kappa)}\int_{\Oint_+^\times \backslash \R_+^n} f(iy) \prod_{j=1}^n y_j^{
  i \rho_j(m)-s + \kappa_j -1/2}dy
\end{align*}
where we have used that the map $y \mapsto 1/y$ maps a fundamental
domain of $\Oint_+^\times \backslash \R_+^n$ to another fundamental
domain. Now (\ref{slffunctional}) follows immediately from the
calculation above since $\sum_{j=1}^n \rho_j(m) = 0$.\qed\newline

Using the Phragm\'{e}n-Lindel\"{o}f principle and the functional equation (\ref{slffunctional}) one obtains that
\begin{align*}
L(1/2 +it,\phi,m) \ll \vert t \vert^{\frac{n}{2}+\epsilon}
\end{align*}
for any $\epsilon > 0$ as $\vert t \vert \to \infty$. This is not enough for our purpose, but any improvement in the exponent will do. In the case $K = \Q$ T. Meurman \cite{meurman90} proved that
\begin{align*}
L(1/2 +it,\phi) \ll \sqrt{r}e^{\pi r/2}\vert t \vert^{\frac{1}{3}+\epsilon},
\end{align*}
where $\frac{1}{4}+r^2$ is the Laplace eigenvalue and the constant implied only depends on $\epsilon$. Recently P. Michel and A. Venkatesh \cite{venkatesh05} and A. Diaconu and P. Garrett \cite{diaconu08} proved the estimate that we need in general:
\begin{theorem}\label{slfsubconv}
There exists some $\delta > 0$ such that
\begin{align*}
L(1/2 +it,\phi,m) \ll \vert t \vert^{\frac{n}{2}-\delta}
\end{align*}
as $\vert t \vert \to \infty$.
\end{theorem}
The generalized Riemann hypothesis implies much more, namely that you can take any $\epsilon > 0$ in the exponent (the Lindel\"{o}f hypothesis for the standard $L$-function). It should be mentioned that the techniques in \cite{petridis01} probably are adequate to provide the subconvexity estimate in Theorem \ref{slfsubconv}.

\section{The Eisenstein Series}
In the case where $K = \Q$ we have the Eisenstein series
\begin{align*}
E(z,s) = \sum_{\gamma \in \Gamma_\infty\backslash \Gamma} \Im (\gamma z)^s.
\end{align*}
In our case of the Hilbert modular group over general $K$ our candidate for the Eisenstein series would be
\begin{equation}\label{preeisenstein}
\sum_{\gamma \in \Gamma_\infty\backslash\Gamma}\prod_{j=1}^n \Im(\gamma^{(j)} z_j)^{s_j}.
\end{equation}
Now for this to be well defined we need every term to be independent of the choice of $\gamma$ in the coset $\Gamma_\infty\backslash\Gamma$. This puts some constraints on the choices of the $s_j$'s. In fact, for (\ref{preeisenstein}) to be well defined it is necessary and sufficient that
\begin{equation}\label{condi}
\vert u^{(1)}\vert^{2s_1} \dots \vert u^{(n)}\vert ^{2s_n} = 1
\end{equation}
for all $u \in \Oint^\times$. The condition (\ref{condi}) is certainly equivalent to
\begin{equation}\label{condiii}
s_1\log \vert \epsilon_j^{(1)}\vert + \dots + s_n \log \vert \epsilon_j^{(n)}\vert = i \pi m_j
\end{equation}
for $j = 1,\dots,n-1$ where $m_j \in \Z$. Let $m =
(m_1,\dots,m_{n-1})\in \Z^{n-1}$ be a fixed vector. If we fix the parameter $s \in \C$ and solve the system of equations
\[
\begin{pmatrix} \log \vert \epsilon_1^{(1)} \vert & \cdots & \log \vert \epsilon_1^{(n)}\vert \\
\cdots & \cdots & \cdots\\
\log \vert \epsilon_{n-1}^{(1)}\vert & \cdots & \log \vert \epsilon_{n-1}^{(n)}\vert \\
1 & \cdots & 1
\end{pmatrix}  \begin{pmatrix} s_1\\ \cdots\\ s_n\end{pmatrix} = \begin{pmatrix} i\pi m_1\\ \cdots\\ i\pi m_{n-1}\\ ns\end{pmatrix}, 
\]
we get the solution (cf. (\ref{unitmatrix}))
\[
s_j = s + i\pi
\sum_{q=1}^{n-1}m_q e_{j,q}=  s+i \rho_j(m)
\]
for $j =1,\dots,n$. From now on we will view $s_j$ as a function of $m$ and $s$. Thus in conclusion we define the Eisenstein series for $\Gamma$ as
\begin{align}
E(z,s,m) = \sum_{\gamma \in \Gamma_\infty\backslash\Gamma}\prod_{j=1}^n \Im(\gamma^{(j)} z_j)^{s_j},
\end{align}
which is absolutely convergent for $\Re(s) > 1$ (cf. \cite{efrat87} p. 42). It was proved in \cite{efrat87} that $E(z,s,m)$ has a meromorphic continuation to the entire $s$-plane, and that $E(z,s,m)$ is holomorphic on the line $\Re(s) = 1/2$.

One can verify that the Eisenstein series is an automorphic form with Laplace eigenvalues $s_j(1-s_j)$ and thus it admits a Fourier expansion. When we calculate the Fourier coefficients it will be convenient to consider the following generalization of the divisor function
\begin{align*}
\sigma_{s,m}(l) = \sum_{\substack{(c) \subset \Oint \\ c \mid l}}\chi_{2m}(c)\vert N(c)\vert^s.
\end{align*}
Note that $\sigma_{s,m}$ only depends on the ideal $(l)$. The Fourier coefficients are known from \cite{efrat87} Section II.2:
\begin{theorem}
For $l \in \Oint$ let $a_l(y,s,m)$ denote the $l$-th Fourier coefficient of
$E(z,s,m)$. For $l \ne 0$ we have that
\begin{align*}
a_l(y,s,m) = \frac{2^n \pi
  ^{ns}\sigma_{1-2s,-m}(l)}{\chi_m(\D)D^s\zeta(2s,-2m)}\prod_{j=1}^n\frac{\sqrt{y_j}K_{s_j-\frac{1}{2}}(2\pi
  \vert (l/\omega)^{(j)}\vert y_j)\vert l^{(j)}\vert^{s_j-\frac{1}{2}}}{\Gamma(s_j)}.
\end{align*}
The zeroth Fourier coefficient is given by
\begin{align*}
a_0(y,s,m) = \left(\prod_{j=1}^n y_j\right)^s\chi_{m}(y)+\phi(s,m)\left(\prod_{j=1}^n y_j\right)^{1-s}\chi_{-m}(y)
\end{align*}
where
\begin{align*}
\phi(s,m) = \frac{\zeta(2s-1,-2m)\pi^{\frac{n}{2}}}{\zeta(2s,-2m)\sqrt{D}}\prod_{j=1}^n\frac{\Gamma(s_j-\frac{1}{2})}{\Gamma(s_j)}.
\end{align*}
Note that $\phi(s,m)$ is unitary for $\Re(s) = \frac{1}{2}$.
\end{theorem}
As in the classical case we also need to consider incomplete
Eisenstein series, i.e. automorphic functions on $\ghn$ formed as Poincar\'{e} series which fail to
be eigenfunctions of the automorphic Laplacian. Let $h \in C_b^\infty(\R_+)$ and assume that $h(y)y^p \to 0$ as $y \to \infty$
and $h(y)y^{-p} \to 0$ as $y \to 0$ for all $p \in \N$. For $m \in \Z^{n-1}$ we define
\begin{align}
F(z, h,m) = \sum_{\gamma \in \Gamma^\infty\backslash
  \Gamma} h\left(\prod_{j=1}^n \Im (\gamma^{(j)}z_j) \right)\prod_{j=1}^n \Im (\gamma^{(j)}z_j)^{i\rho_j(m)}.
\end{align}
We will refer to $F(z,h,m)$ as the
incomplete Eisenstein series induced by $h$ with parameter $m$. One easily checks that the incomplete Eisenstein series decay faster
than any polynomial in the cusp. In particular they are square integrable since
they are bounded. Choosing explicit representatives we see that
\begin{gather}
\begin{split}\label{iesexp}
F(z,h,0) ={} & h\left(\prod_{j=1}^n y_j\right) + h\left(\prod_{j=1}^n \frac{y_j}{x_j^2+y_j^2}\right) +\\
 &\frac{1}{2}\sum_{\substack{c,d \in \Oint^\times\backslash \Oint^*\\ \gcd(c,d) = 1}}h\left(\prod_{j=1}^n\frac{y_j}{(c^{(j)}x_j + d^{(j)})^2 + (c^{(j)}y_j)^2}\right).
\end{split}
\end{gather}
The following proposition reflects the fact that the Hecke $L$-function $\zeta(s,m)$
has a pole at $s=1$ if $m
=0$ but is regular at $s=1$ if $m \ne 0$:
\begin{proposition}\label{intcalc}
For $m \ne 0$ we have
\begin{align*}
\int_{\ghn}F(z,h,m)d\mu(z) = 0.
\end{align*}
We also have
\begin{align*}
\int_{\ghn}F(z,h,0) d\mu(z)= 2^{n-1}R \sqrt{D}\int_0^\infty \frac{h(w)}{w^2}dw.
\end{align*}
\end{proposition}
\proof The last statement follows immediately from change of variables
using the injective map from Proposition \ref{fundprop} and Lemma \ref{jacob}.

The first statement follows from a similar argument. Using again the
map from Proposition \ref{fundprop} and the relation (\ref{fundequa})
we are lead to consider the integral (which only differs from the
integral we wish to compute by scaling with a factor of $R$)
\begin{align*}
&\int_{\R_+\times [-1,1]^{n-1}}
\frac{h(\widetilde{y}_1)}{\widetilde{y}_1^2}\exp\left(i\pi\sum_{q=1}^{n-1}m_q\sum_{i=2}^n
\widetilde{y}_i\sum_{j=1}^n e_{j,q} \log\vert \epsilon_{i-1}^{(j)}\vert\right)d\widetilde{y}\\
&=\int_{\R_+\times [-1,1]^{n-1}}
\frac{h(\widetilde{y}_1)}{\widetilde{y}_1^2}\exp\left(i\pi\sum_{q=1}^{n-1}m_q\widetilde{y}_{q+1}\right)d\widetilde{y}.
\end{align*}
From this the statement is obvious.\qed\newline

The space spanned by incomplete Eisenstein series will be denoted
$\eisen(\ghn)$. Using the transformation from Proposition
\ref{intcalc} it is clear that the orthogonal complement to
$\eisen(\ghn)$ is the set of functions $f \in L^2(\ghn)$ for which
\begin{align}\label{fouriervanish}
\int_F f(z)dx = 0,
\end{align}
i.e. the zeroth Fourier coefficient vanishes. As in the classical
case $K = \Q$ the space $\eisen(\ghn)^\perp$ is the closure of the
space spanned by cusp forms $\cusp(\ghn)$ (see \cite{efrat87} Theorem II.9.8).
Thus we have the decomposition:
\begin{align}\label{specres}
L^2(\ghn) = \overline{\cusp(\ghn)} \oplus
\overline{\eisen(\ghn)}.
\end{align}
Note that the functions in $\cusp(\ghn)$ are orthogonal
to the constant functions.
\section{Quantum Unique Ergodicity}
We wish to investigate the behaviour of the measure
\begin{align*}
d\mu_{m,t} = \vert E(z,1/2 + it,m)\vert^2d\mu
\end{align*}
as $t \to \infty$. This is the large eigenvalue limit, since the Laplace eigenvalue of $E(z,1/2 + it,m)$ is $nt^2 + n/4 +\sum_{j=1}^n \rho_j(m)^2$.

In the subsequent sections we will prove the following two results:
\begin{theorem}\label{QUEEISEN}
Consider an incomplete Eisenstein series $F(z,h,k)$. Then we have that 
\begin{align}\label{maincontrib}
\frac{1}{\log t}\int_{\ghn} F(z,h,k) d\mu_{m,t}(z) \to \frac{\pi^{n} nR}{2 D \zeta(2,0)} \int_{\ghn} F(z,h,k) d\mu(z)
\end{align}
as $t \to \infty$. Note in particular that for $k \ne 0$
\begin{align}
\frac{1}{\log t}\int_{\ghn} F(z,h,k) d\mu_{m,t}(z) \to 0
\end{align}
as $t \to \infty$, cf. Proposition \ref{intcalc}.
\end{theorem}
It is interesting that the asymptotics in (\ref{maincontrib}) do not depend on $m$. The constant $\frac{\pi^{n} nR}{ 2 D\zeta(2,0)}$ can also be given in terms of the volume, since (see \cite{garrett90})
\begin{align}
\mu(\ghn) = \frac{2\zeta(2,0)D^{\frac{3}{2}}}{\pi^n}.
\end{align}

Note that since $\zeta(2) = \frac{\pi^2}{6}$ the result above reduces
to the result found by W. Luo and P. Sarnak in \cite{luo95} for $K = \Q$. The results
differ by a factor of $16$ -- they obtain the asymptotics
\begin{align}\label{maincontribls}
\int_{\ghn} F(z,h) d\mu_{t}(z) \sim \frac{48}{\pi}  \log t \int_{\ghn} F(z,h) d\mu(z)
\end{align}
as $t \to \infty$. This difference is due to a disagreement regarding the value of the integral (\ref{besselint}) below, which exactly accounts for the factor of $16$. In this connection two other errors in \cite{luo95} should be mentioned.
A factor of $2$ is missing in the Fourier expansion of the
Eisenstein series on page 211. This error is cancelled though since a
factor of $\frac{1}{2}$ is missing in front of the logarithmic
derivatives of $\Gamma(s/2 \pm it)$ on page 216.

We also obtain the asymptotics for primitive cusp forms:
\begin{theorem}\label{QUECUSP}
Let $\phi$ be a primitive cusp form. Then
\begin{align}
\int_{\ghn} \phi(z) d\mu_{m,t}(z) \to 0
\end{align}
as $t \to \infty$.
\end{theorem}
Combining Theorem \ref{QUEEISEN} and Theorem \ref{QUECUSP} we can now prove Theorem \ref{quelemma}:
\proof[Proof of Theorem \ref{quelemma}] Let $\epsilon > 0$ be given and set $\Theta = \frac{\pi^{n} nR}{2D\zeta(2,0)}$. One can prove that the
functions which are a sum of a finite number of primitive cusp forms and
incomplete Eisenstein series are dense in the space of
continuous functions which vanish in the cusp $C_0(\ghn)$ equipped
with the sup norm. Hence let $F \in C_c(\ghn)$ and choose primitive cusp forms
$g_1,\dots,g_k$, functions $h_1,\dots,h_l \in C_c^\infty(\R_+)$ and
parameters $m_1,\dots,m_l$ such that
\begin{align*}
\Vert F-G \Vert_\infty \le \frac{\epsilon}{2M\mu(\ghn)},
\end{align*}
where $G(z) = \sum_{j=1}^k g_j(z) + \sum_{i=1}^l F(z,h_i,m_i)$ and $M$ is a constant depending on the field $K$ -- in the case $K = \Q$ one can choose $M = 4$. Now
since cusp forms decay exponentially in the cusp it follows from (\ref{iesexp}) that we can choose a non-negative $h \in
C^\infty(\R_+)$ of sufficiently rapid decay such that
\begin{align*}
\vert F(z)-G(z) \vert \le F(z,h,0)  <
\frac{\epsilon}{2\mu(\ghn)}
\end{align*}
for all $z \in \ghn$. Thus by Theorem \ref{QUEEISEN}
\begin{align*}
\limsup_{t \to \infty} \frac{1}{\Theta \log t} \left\vert \int_\ghn
(F(z)-G(z)) d\mu_{m,t} (z) \right\vert < \frac{\epsilon}{2}.
\end{align*}
Theorem \ref{QUEEISEN} and Theorem \ref{QUECUSP} give us that
\begin{align*}
\lim_{t \to \infty} \frac{1}{\Theta \log t}\int_{\ghn}G(z)d\mu_{m,t}(z) = \int_{\ghn}G(z)d\mu(z).
\end{align*}
Hence
\begin{align}\label{compacteq}
\limsup_{t \to \infty} \left\vert \frac{1}{\Theta \log t}\int_\ghn
F(z) d\mu_{m,t} (z) -\int_\ghn
F(z) d\mu (z)\right\vert < \epsilon.
\end{align}
This proves the theorem, since (\ref{compacteq}) holds for any $\epsilon > 0$.\qed\newline

Finally, this enables us to prove the main theorem:
\proof[Proof of Theorem \ref{que}] Let $F,G,f,g \in C_c(\ghn)$ be chosen such that
\begin{align*}
F \ge 1_{A} \ge f \ge 0 
\end{align*}
and
\begin{align*}
G \ge 1_{B} \ge g \ge 0,
\end{align*}
where $1_A$ denotes the indicator function. Then 
\begin{align*}
\frac{\int_{\ghn}f(z)d\mu_{m,t}(z)}{\int_{\ghn}G(z)d\mu_{m,t}(z)} \le \frac{\mu_{m,t}(A)}{\mu_{m,t}(B)} \le \frac{\int_{\ghn}F(z)d\mu_{m,t}(z)}{\int_{\ghn}g(z)d\mu_{m,t}(z)}.
\end{align*}
By Theorem \ref{quelemma} we see that
\begin{align*}
\frac{\int_{\ghn}f(z)d\mu(z)}{\int_{\ghn}G(z)d\mu(z)} \le \liminf_{t
  \to \infty} \frac{\mu_{m,t}(A)}{\mu_{m,t}(B)} \le \limsup_{t
  \to \infty} \frac{\mu_{m,t}(A)}{\mu_{m,t}(B)} \le \frac{\int_{\ghn}F(z)d\mu(z)}{\int_{\ghn}g(z)d\mu(z)}.
\end{align*}
Since this holds for all $F$, $G$, $f$ and $g$ the result follows.\qed
\section{Proof of Theorem \ref{QUEEISEN}}
Consider $F(z,h,k) \in \eisen(\Gamma \backslash \UH^n)$. By standard unfolding arguments we see that
\begin{align*}
\int_{\Gamma\backslash \UH^n} F(z,h,k)& d\mu_{m,t}\\
&= \int_{\Gamma\backslash
 \UH^n}F(z,h,k)\vert E(z,1/2 + it,m)\vert^2 \frac{dxdy}{y_1^2\dots y_n^2}\\
&= \int_{U_\infty} h\left(\prod_{j=1}^n y_j\right)\int_F
\vert E(z,1/2+it,m)\vert^2 \frac{dxdy}{\prod_{j=1}^n y_j^{2-i\rho_j(k)}}.
\end{align*}
Using the Fourier expansion of the Eisenstein series we get
\begin{align*}
\frac{1}{\sqrt{D}}\int_{F} \vert E(z,1/2+it,m)\vert^2 dx ={}& 2 \prod_{j=1}^n y_j + 2\Re
\left(\prod_{j=1}^n y_j^{1+2it}\chi_{2m}(y)
\overline{\phi\left(1/2+it,m\right)}\right) +\\
&\frac{4^n \pi^n\prod_{j=1}^n y_j}{D\vert
  \zeta(1+2it,-2m)\vert^2}\sum_{l \in
  \Oint^*}\vert\sigma_{-2it,-m}(l)\vert ^2\times\\
&\prod_{j=1}^n\frac{\bigl\vert K_{it+i\rho_j(m)} \left(2\pi\vert (\omega^{-1} l)^{(j)}\vert
  y_j\right)\bigr\vert^2}{\vert\Gamma(1/2+it +i\rho_j(m))\vert^2}.
\end{align*}
Now write
\begin{align*}
\int_{\Gamma\backslash \UH^n} F(z,h,k) d\mu_{m,t} = F_1(t)+F_2(t)
\end{align*}
where
\begin{align*}
F_1&(t) = 2\sqrt{D}\int_{U_\infty} h\left(\prod_{j=1}^n y_j\right)\times\\
&\left(\prod_{j=1}^n y_j + \Re
\left(\prod_{j=1}^n y_j^{1+2it}\chi_{2m}(y)
\overline{\phi\left(1/2+it,m\right)}\right)\right)\frac{dy}{\prod_{j=1}^n
y_j^{2-i\rho_j(k)}}
\end{align*}
and
\begin{align*}
F_2(t) ={}& \frac{4^n \pi^n}{\sqrt{D}\vert
  \zeta(1+2it,-2m)\vert^2}\sum_{l \in
  \Oint^*}\int_{U_\infty} h\left(\prod_{j=1}^n y_j\right)\vert\sigma_{-2it,-m}(l)\vert^2\times\\ 
&\prod_{j=1}^n\frac{\bigl\vert K_{it+i\rho_j(m)} \left(2\pi\vert (\omega^{-1} l)^{(j)}\vert
  y_j\right)\bigr\vert^2}{\vert\Gamma(1/2+it +i\rho_j(m))\vert^2}\frac{dy}{\prod_{j=1}^n
y_j^{1-i\rho_j(k)}}\\
={}&\frac{4^n \pi^n}{\sqrt{D}\vert
  \zeta(1+2it,-2m)\vert^2}\sum_{l \in
  \Oint_+^\times \backslash \Oint^*}\int_{\R_+^n} h\left(\prod_{j=1}^n y_j\right)\vert\sigma_{-2it,-m}(l)\vert ^2\times\\ 
&\prod_{j=1}^n\frac{\bigl\vert K_{it+i\rho_j(m)} \left(2\pi\vert (\omega^{-1} l)^{(j)}\vert
  y_j\right)\bigr\vert^2}{\vert\Gamma(1/2+it
  +i\rho_j(m))\vert^2}\frac{dy}{\prod_{j=1}^n
y_j^{1-i\rho_j(k)}}.
\end{align*}
It is clear that $F_1(t)$ is a bounded function of $t$.

Before we go on we
need to consider a new $L$-function. For $a$ purely imaginary we associate to $\sigma_{a,m}$ an
$L$-function which can be computed in terms of $\zeta(s,m)$:
\begin{align*}
&\sum_{\liea \ne 0} \frac{\chi_{m'}(\liea)\vert \sigma_{a,m}(\liea)
  \vert^2}{\idnorm(\liea)^s} = \prod_{\liep} \sum_{k=0}^\infty \frac{\chi_{m'}(\liep)^k\sigma_{a,m}(\liep^k)\sigma_{-a,-m}(\liep^k)
  }{\idnorm(\liep)^{ks}}\\
&=\prod_{\liep} \sum_{k=0}^\infty
  \frac{\chi_{m'}(\liep)^k}{\idnorm(\liep)^{ks}}\frac{1-\chi_{2m}(\liep)^{k+1}\idnorm(\liep)^{a(k+1)}}{1-\chi_{2m}(\liep)\idnorm(\liep)^a}\frac{1-\chi_{-2m}(\liep)^{k+1}\idnorm(\liep)^{-a(k+1)}}{1-\chi_{-2m}(\liep)\idnorm(\liep)^{-a}}\\
&=
  \prod_{\liep}\frac{1}{(1-\chi_{-2m}(\liep)\idnorm(\liep)^{-a})(1-\chi_{2m}(\liep)\idnorm(\liep)^a)}\times\\
&\phantom{=\vert}\sum_{k=0}^\infty (2 \chi_{m'}(\liep)^k\idnorm(\liep)^{-sk}
  -\chi_{m'}(\liep)^k\chi_{2m}(\liep)^{k+1}\idnorm(\liep)^{(a-s)k+a}-\\
&\phantom{=\vert}\chi_{m'}(\liep)^k\chi_{-2m}(\liep)^{k+1}\idnorm(\liep)^{-(a+s)k-a})\\
&=\prod_{\liep}\frac{1}{(1-\chi_{-2m}(\liep)\idnorm(\liep)^{-a})(1-\chi_{2m}(\liep)\idnorm(\liep)^a)}\times\\
&\phantom{=}\left(\frac{2}{1-\chi_{m'}(\liep)\idnorm(\liep)^{-s}}-\frac{\chi_{2m}(\liep)\idnorm(\liep)^a}{1-\chi_{m'+2m}(\liep)\idnorm(\liep)^{a-s}}-\frac{\chi_{-2m}(\liep)\idnorm(\liep)^{-a}}{1-\chi_{m'-2m}(\liep)\idnorm(\liep)^{-a-s}}
 \right)\\
&=\prod_{\liep}\frac{1+\chi_{m'}(\liep)\idnorm(\liep)^{-s}}{(1-\chi_{m'}(\liep)\idnorm(\liep)^{-s})(1-\chi_{m'+2m}(\liep)\idnorm(\liep)^{a-s})(1-\chi_{m'-2m}(\liep)\idnorm(\liep)^{-a-s})}\\
&=\frac{\zeta(s,m')^2\zeta(s-a,m'+2m)\zeta(s+a,m'-2m)}{\zeta(2s,2m')}.
\end{align*}
To deal with
$F_2(t)$ we consider the Mellin transform $Mh$ of $h$, i.e.
\[
(M h)(r) =  \int_0^\infty h(w) w^{-r-1} dw.
\]
Note that we have the opposite sign convention in the definition of
the Mellin transform than the usual one. However, this is also the convention used in \cite{luo95}, and it is the practical one since we then
avoid considering $\zeta(-s,m)$ on the left half plane. By the Mellin inversion formula we have
\[
h(w) = \frac{1}{2\pi i}\int_{(\sigma)}(M h)(r) w^r dr
\]
for all $\sigma \in \R$. Thus using the $L$-function we considered earlier we can rewrite the integral $F_2(t)$ as
\begin{align*}
F_2(t)\\
={}& \frac{(4\pi)^n}{2\pi i \sqrt{D}\vert
  \zeta(1+2it,-2m)\vert^2}\sum_{l \in
  \Oint_+^\times\backslash\Oint^*}\int_{\R_+^n}\int_{(2)}(M h)(r)\vert\sigma_{-2it,-m}(l)\vert ^2 \times\\ 
&\prod_{j=1}^n\frac{\bigl\vert K_{it+i\rho_j(m)} \left(2\pi\vert (\omega^{-1} l)^{(j)}\vert
  y_j\right)\bigr\vert^2}{\vert\Gamma(1/2+it
  +i\rho_j(m))\vert^2}y_j^{i\rho_j(k)+r-1}drdy\displaybreak\\
={}& \frac{(4\pi)^n}{2\pi i \sqrt{D}\vert
  \zeta(1+2it,-2m)\vert^2 \prod_{j=1}^n\vert\Gamma(1/2+it
  +i\rho_j(m))\vert^2}\int_{(2)}(M
  h)(r)\times\\
&\sum_{l \in
  \Oint_+^\times\backslash\Oint^*}\
\vert\sigma_{-2it,-m}(l)\vert ^2\int_{\R_+^n}\prod_{j=1}^n\bigl\vert K_{it+i\rho_j(m)} \left(2\pi\vert (l/\omega)^{(j)}\vert
  y_j\right)\bigr\vert^2y_j^{i\rho_j(k)+r-1}dydr\\
={}&\frac{(4 \pi)^n}{2 \pi i 2^{3n} \sqrt{D}\vert
  \zeta(1+2it,-2m)\vert^2 \prod_{j=1}^n\vert\Gamma(1/2+it
  +i\rho_j(m))\vert^2}\int_{(2)}(M
  h)(r)\times\\
&\sum_{l \in
  \Oint_+^\times\backslash\Oint^*}\vert\sigma_{-2it,-m}(l)\vert ^2\prod_{j=1}^n
  \frac{\vert \omega^{(j)} \vert^{i\rho_j(k)+r}
  \Gamma((i\rho_j(k)+r)/2)^2}{
  \pi^{i\rho_j(k)+r} \vert l^{(j)} \vert^{i\rho_j(k)+r}\Gamma(i\rho_j(k)+r)}\times \\
&\Gamma((i\rho_j(k)+r)/2+it+i\rho_j(m)) \Gamma((i\rho_j(k)+r)/2-it-i\rho_j(m)) dr\\
={}&\frac{(4 \pi)^n}{2\pi i 2^{3n} \sqrt{D}\vert
  \zeta(1+2it,-2m)\vert^2 \prod_{j=1}^n\vert\Gamma(1/2+it
  +i\rho_j(m))\vert^2}\int_{(2)}B_k(r,t,h)dr
\end{align*}
where
\begin{align*}
B_k(r,t,h)={}&(M h)(r)
  \frac{\zeta(r,-k)^2\zeta(r+2it,-k-2m)\zeta(r-2it,-k+2m)}{\zeta(2r,-2k)\pi^{nr}}\times\\
&\prod_{j=1}^n
  \frac{\vert \omega^{(j)} \vert^{i\rho_j(k)+r}\Gamma((i\rho_j(k)+r)/2)^2}{\Gamma(i\rho_j(k)+r)}\times\\
&\Gamma((i\rho_j(k)+r)/2+it+i\rho_j(m)) \Gamma((r+i\rho_j(k))/2-it-i\rho_j(m)).
\end{align*}
Note that we have used the fact that for any $b \in \R$ we
have the formula (see \cite{iwaniec02} Section B.4)
\begin{equation}
\int_0^\infty \vert K_{ib}(2\pi t) \vert^2 t^{s-1} dt =\frac{
  \Gamma(s/2+ib) \Gamma(s/2-ib) \Gamma(s/2)^2}{2^3 \pi^s \Gamma(s)}.
\end{equation}
Clearly $(Mh)(r)$ is bounded for $\frac{1}{2} \le \Re(r) \le 2$ and
$\Gamma$ decays exponentially in vertical strips by Stirling's formula. Furthermore $\zeta(\sigma + it,k)$ is
polynomially bounded in $t$ for $\frac{1}{2} \le \sigma \le 2$. Hence
we can move the integration from the vertical line $\Re(r) = 2$ to
the vertical line $\Re(r) = \frac{1}{2}$ perhaps picking up residues from poles at $r = 1$ and $r=1\pm 2it$:
\begin{align*}
F_2(t)={}&\frac{(\pi/2)^n\int_{(1/2)}B_k(r,t,h)dr}{2\pi i \sqrt{D}\vert
  \zeta(1+2it,-2m)\vert^2 \prod_{j=1}^n\vert\Gamma(1/2+it
  +i\rho_j(m))\vert^2}+\\
&\frac{(\pi/2)^n\res_{r = 1}B_k(r,t,h)}{\sqrt{D}\vert
  \zeta(1+2it,-2m)\vert^2 \prod_{j=1}^n\vert\Gamma(1/2+it
  +i\rho_j(m))\vert^2}+O(t^{-10})
\end{align*}
where the $O(t^{-10})$ term comes from the possible residues from poles at
$r=1\pm 2it$, since $(Mh)(\sigma + it)$ is of rapid decay as $t \to
\infty$. Let us evaluate the first term. Since Stirling's formula is no good near the real axis in our case, we have to work around that. Note that for $a,b \in \R$ we have
\begin{align*}
e^{-\vert a+b\vert}e^{-\vert a-b\vert} \le e^{-2\vert a\vert}.
\end{align*}
If $\vert a + b \vert \ge 1$ and $a \ne 0$ we also have that
\begin{align*}
\frac{1}{\vert a + b\vert} \le \frac{1 + \vert b \vert}{\vert a \vert}.
\end{align*}
We can now evaluate the first term. Since we are only interested in the asymptotics as $t \to \infty$ we can assume that $t \ge 1$. Using the subconvexity estimate from
Theorem \ref{sohneest} and Stirling's formula we see that ($C_1,C_2,C_3 > 0$ are suitable constants)
\begin{align*}
\int_{(1/2)}\vert B_k & (r,t,h) \vert dr \le\\
& e^{-\pi t n}t^{-\frac{n}{6}+\epsilon}C_1 \int_{-\infty}^\infty \vert (Mh)(1/2 + i w)\vert \left(1+\vert w\vert\right)^{\frac{2n}{3}+\epsilon}dw+\\
& e^{-\pi t n}t^{-\frac{n}{4}+\epsilon}C_2 \int_{2(t+\rho_j(m)-1)-\rho_j(k)}^{2(t+\rho_j(m)+1)-\rho_j(k)} \vert (Mh)(1/2 + i w)\vert dw+\\
& e^{-\pi t n}t^{-\frac{n}{4}+\epsilon}C_3 \int_{-2(t+\rho_j(m)+1)-\rho_j(k)}^{-2(t+\rho_j(m)-1)-\rho_j(k)}\vert (Mh)(1/2 + i w)\vert dw.
\end{align*}
Since $Mh$ is of rapid decay the first term dominates, and we obtain the estimate
\begin{align*}
\int_{(1/2)}B_k(r,t,h)dr \ll e^{-t\pi n}\vert t\vert^{-\frac{n}{6}+\epsilon}.
\end{align*}
By Corollary \ref{convest2} and Stirling's formula we see that
\begin{align*}
\frac{\int_{(1/2)}B_k(r,t,h)dr}{\vert
  \zeta(1+2it,-2m)\vert^2 \prod_{j=1}^n\vert\Gamma(1/2+it
  +i\rho_j(m))\vert^2} \ll \vert t
  \vert^{-\frac{n}{6}+\epsilon} 
\end{align*}
for any $\epsilon > 0$.

Now we turn to the residue term. Since $\zeta(s,k)$ is regular at $s =
1$ for $k \ne 0$ the resiude term will vanish in this case and we are done. Assume
therefore that $k = 0$. We know that
\begin{align*}
\zeta(s,0) = \frac{\zeta_{-1}}{s-1} + \zeta_0 + O(s-1)
\end{align*}
and hence
\begin{align*}
\zeta(s,0)^2 = \frac{\zeta_{-1}^2}{(s-1)^2} +
\frac{2\zeta_{-1}\zeta_0}{s-1}+ O(1)
\end{align*}
as $s \to 1$ where $\zeta_{-1} = \frac{2^{n-1}R}{\sqrt{D}}$ and $\zeta_0$ is some constant. Now introduce $G(r,t,h)$ defined by
\begin{align*}
B_0(r,t,h) = \zeta(r,0)^2G(r,t,h).
\end{align*}
We see that
\begin{align*}
\res_{r=1} B_0(r,t,h) = G(1,t,h)\zeta_{-1}\left(2\zeta_0 + \zeta_{-1}\frac{G'(1,t,h)}{G(1,t,h)}\right).
\end{align*}
Note that
\begin{align*}
G(1,t,h) = \frac{(Mh)(1)\vert
  \zeta(1-2it,2m)\vert^2}{\zeta(2,0)\pi^n}D\Gamma(1/2)^{2n}\prod_{j=1}^n\vert\Gamma(1/2+it+i\rho_j(m))\vert^{2}
\end{align*}
and
\begin{align*}
\frac{G'(1,t,h)}{G(1,t,h)} ={}&\frac{\zeta'(1+2it,-2m)}{\zeta(1+2it,-2m)}+\frac{\zeta'(1-2it,2m)}{\zeta(1-2it,2m)}+\\
&\frac{1}{2}\sum_{j=1}^n\left(\frac{\Gamma'(1/2+it+i\rho_j(m))}{\Gamma(1/2+it+i\rho_j(m))}+\frac{\Gamma'(1/2-it-i\rho_j(m))}{\Gamma(1/2-it-i\rho_j(m))}
\right) + C
\end{align*}
where $C$ is a constant that does not depend on $t$. Since
\begin{align*}
(Mh)(1) = \frac{2^{1-n}}{\sqrt{D}R} \int_\ghn F(z,h,0)d\mu(z)
\end{align*}
by Proposition \ref{intcalc} we
see using Corollary \ref{colemancor} and Stirling's formula that
\begin{align*}
\frac{1}{\log t}F_2(t) \to \frac{\pi^{n} nR}{2 D \zeta(2,0)} \int_\ghn
F(z,h,0)d\mu(z)
\end{align*}
as $t \to \infty$.
\section{Proof of Theorem \ref{QUECUSP}}
Let $\phi$ be a primitive cusp form with eigenvalues $\frac{1}{4}+r_j^2$ of the
Laplacians $\Delta_j$ and Hecke eigenvalues $\lambda(\liea)$.

We wish to investigate the asymptotic behaviour of the integral
\begin{equation}
\int_\ghn \phi(z)d\mu_{m,t} = \int_\ghn \phi(z) E(z,1/2+it,m)E(z,1/2-it,-m)d\mu
\end{equation}
where we have used the fact that $\overline{E(z,s,m)} =
E(z,\overline{s},-m)$. To this end we consider the integral
\begin{equation}
I(s) = \int_\ghn \phi(z)E(z,1/2+it,m)E(z,s,-m) d\mu
\end{equation}
for $\Re(s) > 1$. We unfold the integral and get using the Fourier
expansions of cusp forms and Eisenstein series that
\begin{align*}
I(s) ={}& \int_{F_\infty} \phi(z)E(z,1/2+it,m)\prod_{j=1}^n y_j^{s_j(-m)-2} dxdy\\
={}& \frac{2^n \pi^{n(1/2+it)}}{\zeta(1+2it,-2m)\chi_m(\D)D^{it}}\int_{U_\infty}\sum_{l \in O^*}\sigma_{-2it,-m}(l)c_l\prod_{j=1}^ny_j^{s_j(-m)-1}\vert l^{(j)}\vert^{it+i\rho_j(m)}\times\\
& \frac{K_{it+i\rho_j(m)}(2\pi \vert (l/\omega)^{(j)}\vert y_j)K_{ir_j}(2\pi \vert (l/\omega)^{(j)}\vert y_j)}{\Gamma(1/2+it+i\rho_j(m))}  dy\\
={}&\frac{2^n
  \pi^{n(1/2+it)}(\prod_{j=1}^n\vert\omega^{(j)}
  \vert^{s-i\rho_j(m)})}{\zeta(1+2it,-2m)\chi_m(\D)D^{it}\prod_{j=1}^n\Gamma(1/2+it+i\rho_j(m))}\sum_{l\in \Oint_+^\times\backslash \Oint^*}\chi_{2m}(l)\times\\
&\idnorm((l))^{it-s}\sigma_{-2it,-m}(l)c_l\int_{\R_+^n}\prod_{j=1}^nK_{it+i\rho_j(m)}(2\pi
y_j)K_{ir_j}(2\pi y_j) y_j^{s_j(-m)-1}dy.
\end{align*}
For $a,b \in \R$ consider the meromorphic function on $\C$:
\begin{align*}
\Gamma &(s,a,b) =\\
&\frac{\Gamma((s+ia+ib)/2)\Gamma((s+ia-ib)/2)\Gamma((s-ia-ib)/2)\Gamma((s-ia+ib)/2)}{2^3\pi^s\Gamma(s)}.
\end{align*}
It is well known (see \cite{iwaniec02} Section B.4) that
\begin{align}\label{besselint}
\int_0^\infty K_{ia}(2\pi
t)K_{ib}(2\pi t) t^{s-1}dt = \Gamma(s,a,b).
\end{align}
So we get
\begin{align*}
I(s) ={}& \frac{2^{n}
  \pi^{n(1/2+it)}}{\zeta(1+2it,-2m)\chi_m(\D)D^{it}}\prod_{j=1}^n\frac{\vert\omega^{(j)} \vert^{s-i\rho_j(m)}\Gamma(s_j(-m),r_j,t+\rho_j(m))}{\Gamma(1/2+it+i\rho_j(m))}\times \\
&R(s)\sum_{\beta \in \Oint_+^\times\backslash \Oint^\times}c_\beta
\end{align*}
where
\begin{align*}
R(s) ={}& \sum_{\liea \subset
  \Oint}\chi_{2m}(\liea)\idnorm(\liea)^{it-s}\sigma_{-2it,-m}(\liea)\lambda(\liea)\\
={}& \prod_{\liep}\sum_{k =
  0}^\infty\chi_{2m}(\liep)^k\idnorm(\liep)^{k(it-s)}\sigma_{-2it,-m}(\liep^k)\lambda(\liep^k)\\
={}& \prod_{\liep}\sum_{k =
  0}^\infty\chi_{2m}(\liep)^k\idnorm(\liep)^{k(it-s)}\lambda(\liep^k)\sum_{j
  = 0}^k\chi_{-2m}(\liep)^j\idnorm(\liep)^{-2ijt}\\
={}& \prod_{\liep}\sum_{k =
  0}^\infty\chi_{2m}(\liep)^k\idnorm(\liep)^{k(it-s)}\lambda(\liep^k)\frac{1-\chi_{-2m}(\liep)^{k+1}\idnorm(\liep)^{-2(k+1)it}}{1-\chi_{-2m}(\liep)\idnorm(\liep)^{-2it}}\\
={}&
  \prod_{\liep}\frac{1}{1-\chi_{-2m}(\liep)\idnorm(\liep)^{-2it}}\biggl(\sum_{k=0}^\infty\chi_{2m}(\liep)^k\idnorm(\liep)^{k(it-s)}\lambda(\liep^k)-\\
&\chi_{-2m}(\liep)\idnorm(\liep)^{-2it}\sum_{k=0}^\infty\lambda(\liep^k)\idnorm(\liep)^{k(-it-s)}\biggr)\\
={}&
  \prod_{\liep}\frac{1}{1-\chi_{-2m}(\liep)\idnorm(\liep)^{-2it}}\times\\
&\biggl(\frac{1}{1-\lambda(\liep)\chi_{2m}(\liep)\idnorm(\liep)^{it-s}+\chi_{2m}(\liep)^2\idnorm(\liep)^{2(it-s)}}-\\
&\frac{\chi_{-2m}(\liep)\idnorm(\liep)^{-2it}}{1-\lambda(\liep)\idnorm(\liep)^{-it-s} + \idnorm(\liep)^{2(-it-s)}}\biggr)\\
={}&\prod_{\liep}\frac{1-\chi_{2m}(\liep)\idnorm(\liep)^{-2s}}{(1-\chi_{2m}(\liep)\lambda(\liep)\idnorm(\liep)^{it-s}+\chi_{2m}(\liep)^2\idnorm(\liep)^{2(it-s)})}\times\\
&\frac{1}{(1-\lambda(\liep)\idnorm(\liep)^{-it-s} + \idnorm(\liep)^{2(-it-s)})}\\
={}&\frac{L(s-it,\phi,2m)L(s+it,\phi,0)}{\zeta(2s,2m)}.
\end{align*}
From this we see that $I(s)$ has an analytic continuation to the
entire $s$-plane, and we wish to investigate the asymptotic behaviour
of $I(1/2-it)$ as $t \to \infty$. From Stirling's formula
we deduce that
\begin{align*}
\prod_{j=1}^n\frac{\Gamma(1/2-it-i\rho_j(m),r_j,t+\rho_j(m))}{\Gamma(1/2+it+i\rho_j(m))}\ll \vert t \vert^{-n/2}
\end{align*}
as $t \to \infty$. Using Proposition \ref{convest2} the proof of
Theorem \ref{QUECUSP} boils down to proving a subconvexity estimate
for $L(s,\phi,2m)$ on the line $\Re(s) = \frac{1}{2}$. More precisely we
need the estimate
\begin{align*}
L(1/2 +it,\phi,2m) \ll \vert t \vert^{\frac{n}{2}-\delta}
\end{align*}
as $\vert t\vert \to \infty$ for some $\delta > 0$, and this follows from Theorem \ref{slfsubconv}. Note that if $\phi$ is odd then $I(1/2 - it) = 0$, since $L(1/2,\phi,0)=0$ by the functional equation.\newline

\end{document}